\documentclass[12pt,a4paper]{article}

\usepackage{graphicx,hyperref,enumitem}
\usepackage{amsmath,amsthm,amssymb,color,bbm}
\usepackage{tocloft,accents,array}

\def\lognd{\log\dfrac{2n}{\dmax}}
\def\srho{\rho_{\mathrm{small}}}
\def\brho{\rho_{0}}
\def\dmin{d_{\mathrm{min}}}
\def\dmax{d_{\mathrm{max}}}

\def\pr{\mathbb P}

\normalsize

\date{\today}
\def\nperp{n_{\scriptscriptstyle\perp}}
\def\yvec{\boldsymbol{y}}
\def\Rmodpi{\Reals/\mkern-1mu\pi}

\def\semiabs#1{\abs{#1}_\pi}

\def\norm#1{\mathopen\|#1\mathclose\|}
\def\onenorm#1{\norm{#1}_1}

\def\infnorm#1{\norm{#1}_\infty}

\def\sp{\operatorname{span}}
\def\({\bigl(}
\def\){\bigr)}
\def\st{:}

\def\ms#1{\mathopen{\{\!\{}#1\mathclose{\}\!\}}}
\def\Ms#1{\mathopen{\bigl\{\!\bigl\{}#1\mathclose{\bigr\}\!\bigr\}}}
\let\leq=\leqslant
\let\geq=\geqslant
\let\le=\leqslant
\let\ge=\geqslant
\newcommand{\medtilde}{\protect\accentset{\sim}}

\newcommand{\medbar}[1]{\overset{\rule{0.5em}{0.4pt}}{#1}}
\newcommand{\Medbar}[1]{\overset{\kern0.1em\rule{0.5em}{0.4pt}}{#1}}

\def\RT{\operatorname{RT}}
\def\ED{\operatorname{ED}}
\def\EO{\operatorname{EO}}
\def\EOG{\operatorname{EOG}}

\def\Integers{\mathbb{Z}}

\newcommand{\stirlingii}{\genfrac{\{}{\}}{0pt}{}}

\def\thetavec{{\boldsymbol{\theta}}}
\def\phivec{{\boldsymbol{\phi}}}
\def\thetapvec{{\boldsymbol{\theta'}}}

\newtheorem{thm}{Theorem}[section]
\newtheorem{cor}[thm]{Corollary}

\newtheorem{lemma}[thm]{Lemma}

\theoremstyle{definition}

\numberwithin{equation}{section}
\allowdisplaybreaks[2]

\let\originalleft\left
\let\originalright\right
\renewcommand{\left}{\mathopen{}\mathclose\bgroup\originalleft}
\renewcommand{\right}{\aftergroup\egroup\originalright}

\def\Abs#1{\bigl\lvert#1\bigr\rvert} \let\Card=\Abs
\def\nfrac#1#2{{\textstyle\frac{#1}{#2}}}
\def\dfrac#1#2{\lower0.12ex\hbox{\large$\textstyle\frac{#1}{#2}$}}
\def\Dfrac#1#2{\raise0.05ex\hbox{\small$\displaystyle\frac{#1}{#2}$}}

\let\eps=\varepsilon

\def\nperp{n_{\scriptscriptstyle\perp}}
\def\trans{^{\mathrm{T}}}

\def\Ex{\operatorname{\mathbf{E}}}

\def\Cov{\operatorname{Cov}}

\def\Reals{{\mathbb{R}}}
\def\Complexes{{\mathbb{C}}}

\def\Naturals{{\mathbb{N}}}

\def\S{\boldsymbol{S}}

\def\X{\boldsymbol{X}}
\def\Y{\boldsymbol{Y}}
\def\Z{\boldsymbol{Z}}
\def\x{\boldsymbol{x}}
\def\y{\boldsymbol{y}}
\def\calF{\mathcal{F}}

\def\thetavec{\boldsymbol{\theta}}
\def\xvec{\boldsymbol{x}}
\def\yvec{\boldsymbol{y}}

\newcommand{\E}[1]{\mathbf E\left[#1\right]}
\newcommand{\EE}[1]{\mathbf E\bigl[#1\bigr]}
\newcommand{\EEE}[1]{\mathbf E\biggl[#1\biggr]}

\def\nicebreak{\vskip 0pt plus 50pt\penalty-300\vskip 0pt plus -50pt }

\usepackage{bm}
\usepackage{float}
\newtheorem{df}[thm]{Definition}

\renewcommand{\P}[1]{\pr\left(#1\right)}

\newcommand{\Exp}[1]{\exp\left(#1\right)}

\renewcommand{\O}[1]{O\left(#1\right)}

\renewcommand{\th}{\theta}

\newcommand{\ka}{\kappa}
\newcommand{\Ka}[1]{\ka\left(#1\right)}

\def\abs#1{\lvert#1\rvert} \let\card=\abs
\def\Abs#1{\bigl\lvert#1\bigr\rvert} \let\Card=\Abs

\newcommand{\bR}[1]{\left(#1\right)}

\newcommand{\bRg}[1]{\biggl(#1\biggr)}
\newcommand{\bRgg}[1]{\Biggl(#1\Biggr)}

\newcommand{\ST}{\tau}

\makeatletter
\def\and{%
  \end{tabular}%
  \begin{tabular}[t]{c}}
\makeatother

\begin{document}

\title{\Large{\bf
Cumulant expansion for counting Eulerian orientations}}
\author{
\hspace{-.3cm}Mikhail Isaev\thanks{Supported by Australian Research Council grant DP250101611.} \\
\hspace{-.3cm}\small School of Mathematics and Statistics\\[-0.3ex]
\hspace{-.3cm}\small University of New South Wales\\[-0.3ex]
\hspace{-.3cm}\small Sydney, NSW, Australia\\[-0.3ex]
\hspace{-.3cm}\small\tt  m.isaev@unsw.edu.au 
\and
\hspace{-.5cm}Brendan D. McKay\footnotemark[1]\\
\hspace{-.5cm}\small School of Computing \\[-0.3ex]
\hspace{-.5cm}\small Australian National University  \\[-0.3ex]
\hspace{-.5cm}\small Canberra ACT 2601, Australia \\[-0.3ex]
\hspace{-.5cm}\small\tt brendan.mckay@anu.edu.au
\bigskip
\and
\hspace{-.5cm}Rui-Ray Zhang \\
\hspace{-.5cm}\small Barcelona School of Economics \\[-0.5ex]
\hspace{-.5cm}\small 08005 Barcelona, Spain  \\[-0.5ex]
\hspace{-.5cm}\small\tt rui.zhang@bse.eu
}

\date{}

\maketitle
\begin{abstract}
An Eulerian orientation is an orientation of the edges of a graph such that every vertex is balanced: 
its in-degree equals its out-degree.
Counting Eulerian orientations corresponds to the crucial partition function in so-called ``ice-type models'' in statistical physics and is known to be hard for general graphs.
For all graphs with good expansion properties and degrees larger than $\log^{8} n$,
we derive an asymptotic expansion for this count that approximates it to precision $O(n^{-c})$
for arbitrary large $c$, where $n$ is the number of vertices.
 The proof relies on a new tail bound for the cumulant expansion of the Laplace transform, which is of independent interest.
\end{abstract}


\section{Introduction}


All the graphs in this paper are simple.
An \emph{Eulerian orientation} of an undirected graph 
is an assignment of directions to all the edges such that, for any vertex,
the number of edges directed inwards equals the number of edges directed outwards. 
Let $\EO(G)$ denote the number of Eulerian orientations of graph $G$ with vertex set $[n]:=\{1,\ldots,n\}$. 
It is well known that $\EO(G)>0$ if and only if 
$G$ is \emph{Eulerian},
that is, all the degrees of $G$ are even.


Schrijver \cite{schrijver1983bounds} showed that   
 \begin{equation}\label{eq:Sch}
        \Dfrac{1}{2^{|E(G)|}} \prod_{i \in [n]}\binom{d_i}{d_i/2}   \leq    \EO(G) \leq \prod_{i \in [n]}\binom{d_i}{d_i/2}^{\!1/2},
 \end{equation}
 where $d_1,\ldots,d_n$ are the degrees of the vertices,
 and $E(G)$ denotes the edge set.
The lower bound in \eqref{eq:Sch}
coincides with 
a heuristic estimate by Pauling \cite{pauling1935structure} around 90 years ago, which is obtained
by ignoring the dependencies among the events that vertices are unbalanced in a random orientation of the edges.
For the regular case when $d_i =d$ for all $i \in [n]$, an improvement to these bounds is given by Las Vergnas~\cite{LasVergnas1990}.
In a recent paper~\cite{IMZentropy}, 
we study  Eulerian orientations of random graphs, and also
show that the
lower bound in~\eqref{eq:Sch} is low by at most a factor of $e^{22n/7}$
for all simple graphs.

There are no known efficient algorithms for finding $\EO(G)$ exactly. 
Schrijver \cite{schrijver1983bounds} 
pointed out that
the computation of $\EO(G)$ can be reduced to the evaluation of the permanent of a certain matrix associated with the graph.
Mihail and Winkler \cite{mihail1996number} 
showed that  computing the number of Eulerian orientations is 
 \#P-complete. They also  
proposed an approximation algorithm based on the Markov chain Monte Carlo 
method for computing the permanent.

Wide interest  in counting Eulerian orientations 
is due to the equivalence to  the partition function of so-called
ice-type models in statistical physics; see, for example,  \cite{welsh1990computational} and references therein. 
Lieb \cite{lieb1967residual} and Baxter \cite{baxter1969f}  derived asymptotic expressions for 
the number of Eulerian orientations of the square and
triangular lattice by the transfer matrix method 
\cite{lieb1967residual, baxter1969f} using eigenvalues of matrices.
No extension of this to higher dimensions is known,
to the best of our knowledge.

A \emph{regular tournament} on $n$ vertices is an Eulerian orientation of the complete graph $K_{n}$. Clearly, $n$ must be odd for a regular tournament to exist. 
McKay \cite{mckay1990asymptotic} established that, for odd $n \rightarrow \infty$,
\[
    \EO(K_n) \sim 
\left(\Dfrac{n}{e}\right)^{1/2}
\left(\Dfrac{2^{n+1}}{\pi n}\right)^{(n-1)/2}.
\]
Adopting McKay's approach to dense graphs with strong mixing properties, Isaev and Isaeva~\cite{isaev2013asymptotic} 
obtained
\[
\EO(G) 
 \sim \widehat{\EO}(G) :=
\Dfrac{ 2^{ |E(G)| } }{ \sqrt{ \ST(G) } }
\left( \Dfrac{2}{\pi} \right)^{(n-1)/2}
\Exp{ -\dfrac14 \sum_{ jk \in G } 
\biggl( \Dfrac{1}{d_{j}} + \Dfrac{1}{d_{k}} \biggr)^{\!2\,}
},
\]
where $\ST(G)$ denotes the number of spanning trees of $G$,
and the summation is over all (unordered) edges $\{j, k\} \in E(G)$.
More precisely, graphs  considered in \cite{isaev2013asymptotic}   
satisfy $h(G)\geq \gamma n$ for some fixed $\gamma>0$, where $h(G)$
is the \textit{Cheeger constant} (also known as isoperimetric number) of~$G$,
defined by
\[
	h(G):=\min
	\biggl\{ \Dfrac{|\partial_G \,U|}{|U|}  \st  U\subseteq V(G), 1\le |U|\le\dfrac n2 \biggr\},
\]
where $\partial_G \,U$ is the set of edges of $G$ with one end in $U$ 
and the other end in $V(G) \setminus U$.  
Isaev, Iyer and McKay~\cite{isaev2020asymptotic} obtained an asymptotic expression
for $\EO(G)$ under the conditions that $h(G) \ge \gamma \dmax$ and $\dmax \ge n^{1/3 + \eps}$
for constants $\gamma,\eps>0$, where $\dmax$ is the maximum degree of $G$. As a corollary, they established
\begin{equation}\label{EO-ext}
\EO(G) = 
\widehat{\EO}(G)
\Exp{ \O{ \Dfrac{n}{\dmax^{2}} \log\Dfrac{2n}{\dmax} } }.
\end{equation}

Note that if $\dmax \ll n^{1/2}$, then the error term 
$\dfrac{n}{\dmax^{2}} \log\dfrac{2n}{\dmax}$ in \eqref{EO-ext}  grows as $n\rightarrow \infty$.
In this paper, we obtain an asymptotic expansion series that gives the value of  $\EO(G)$ up to precision $O(n^{-c})$ 
for arbitrary large constant $c$ under the assumptions 
that $h(G) \geq \gamma \dmax$ and  $\dmax\gg \log^{8}n$. 
As a corollary, we prove that bound \eqref{EO-ext} also holds for such graphs.  

\subsection{Main results}

We will estimate $\EO(G)$ using a cumulant expansion with respect to a certain Gaussian random vector associated
with $G$. 
We first recall the definition of cumulants.  Let $X_1,\ldots,X_r$ be random variables
with finite moments,
the  {\em joint cumulant} (or {\em mixed cumulant}) is defined by 
\begin{equation}
    \kappa (X_1,\ldots,X_r) := [t_1 \cdots t_r] 
    \log\biggl( \EEE{ \exp\bRg{ \sum_{i=1}^{r} t_i X_i } } \biggr),
    \label{cumu}
\end{equation}
where $[t_1 \cdots t_r]$ stands for the coefficient of $t_1 \cdots t_r$ in the formal series expansion.
If all the random variables $X_1,\ldots,X_r$ are the same variable $X$,
then we write 
\[
   \kappa_r(X):=\kappa(\underbrace{X,\ldots,X}_{r \text{ times}}),
\]
which is the {\em cumulant} of order $r$ for random variable $X$.

Given a graph $G$,  the 
{\it Laplacian matrix} $L=L(G)$ is  the symmetric matrix defined by
\begin{equation}\label{Ldef}
     \xvec\trans\! L\xvec 
     =  \sum_{jk\in G}  (x_j-x_k)^2,
\end{equation}
for $\xvec=(x_1,\ldots,x_n)\trans\in\mathbb R^n$.
 Clearly, $L$ is positive-semidefinite.  
 The  vector 
 $\boldsymbol{1}=(1,\ldots,1)\trans $ is an eigenvector of $L$ with eigenvalue~$0$. If $G$ is connected, then all other eigenvalues of $L$ are positive.  In this case, 
let  $\X_G$ denote  $n$-dimensional singular Gaussian  random vector on
 the subspace 
 \begin{equation}\label{eq:Vn}
    \mathcal{V}_n:= \{\xvec\in \Reals^n \st x_1 +\cdots +x_n = 0\}
 \end{equation}
with density proportional to 
$\Exp{-\frac12\xvec\trans\! L\xvec}$.

Let $( c_{2t} )_{t \ge 1}$  denote the coefficients of the Taylor expansion of $\log\cos x$ at $x=0$, i.e.,
\[
    \log \cos x 
    = \sum_{t\geq 1} c_{2t} x^{2t}
    = - \dfrac12 x^2 - \dfrac{1}{12} x^4 - \dfrac{1}{45} x^6 -\dfrac{17}{2520}x^8 - \cdots,
\]
where
\begin{equation}
    \label{eq:cDef}
  c_{2t} 
  = -
  \frac{4^t (4^{t}-1)
  | B_{2t} |}{2t (2t)!},    
\end{equation}
with 
$B_{2t}$ denoting the Bernoulli number.
For an integer $K \ge 2$, define polynomial $f_{K}$ by
\begin{equation}
f_{K}(\x)
:= \sum_{t=2}^K c_{2t} \sum_{jk \in G} (x_j-x_k)^{2t}.
\label{eq:eoTaylor}
\end{equation}

Our main result is the following theorem.
\begin{thm}
\label{tm-eo}
Let $G = G(n)$ be a graph with even degrees.
Assume that  
\begin{itemize}\itemsep=0pt
\item[(i)]  $\dmax \gg \log^{8} n$, where $\dmax$ is the maximum degree  of $G$;
\item[(ii)] the Cheeger constant $h(G) \ge \gamma \dmax$ for some constant $\gamma > 0$.
\end{itemize}
Let  $c > -1$ be a constant, and 
let $M$ and $K$ be defined by
\[
M := \biggl\lfloor \Dfrac{(c+1)\log n}{\log \dmax - 2\log\log n} \biggr\rfloor
\quad\text{and}\quad
K := \biggl\lceil \Dfrac{(c+1)\log n}{\log \dmax - \log\log n} \biggr\rceil.
\]
Then
as $n \to \infty$,
\begin{equation}
\EO(G) = 
\Dfrac{ 2^{ |E(G)| } }{ \sqrt{ \ST(G) } }
\left( \Dfrac{2}{\pi} \right)^{(n-1)/2}
\exp\bRgg{ 
\sum_{ r =1 }^M \Dfrac{1}{r!}\, \kappa_r \left( f_{K}(\X_G) \right)
+ \O{ n^{-c} } },
\label{eo}
\end{equation}
where 
$\ST(G)$ denotes the number of spanning trees of graph $G$,
and $|E(G)|$ denotes the number of edges in graph $G$.
\end{thm}

As a corollary of Theorem \ref{tm-eo}, we obtain the following result.

\begin{cor}\label{cor-eo}
The asymptotic bound of \eqref{EO-ext} holds under the assumptions of Theorem \ref{tm-eo}. 
\end{cor}
The proofs of  Theorem \ref{tm-eo} 
and Corollary \ref{cor-eo}
are given in Section \ref{S:proof-main}.  
We believe that the assumption that $\dmax\gg \log^{8} n$ is an artefact of our proof techniques and the result can be extended 
to the range $\dmax\gg \log^2n$ at least.

\subsection{Overview of the paper}
\label{S:overview}

To help with the navigation within the paper, we briefly present our approach here which builds on previous works \cite{mckay1990asymptotic,isaev2013asymptotic,isaev2020asymptotic}.
The number of Eulerian orientations of a graph $G$ is the 
constant term in
\[
    \prod_{jk\in G} \,\Bigl( \Dfrac{x_j}{x_k}+\Dfrac{x_k}{x_j}\Bigr),
\]
where the product is over the edges $\{j,k\}$ of~$G$.
Applying Cauchy's theorem with the substitution $x_j=e^{i\theta_j}$,
we obtain
\begin{equation}\label{EOinitial1}
   \EO(G) = \Dfrac{2^{\abs{E(G)}}}{(2\pi)^n}
      \int_{-\pi}^\pi\!\cdots\!\int_{-\pi}^\pi 
      \prod_{jk\in G} \cos(\theta_j-\theta_k)\,d\thetavec,
\end{equation}
where $\thetavec=(\theta_1,\ldots,\theta_n)$.
If $G$ is Eulerian, the identity $\cos(\pi+\theta)=-\cos \theta$
preserves the integrand, so we can simply~\eqref{EOinitial1} to
\begin{equation}\label{EOinitial2}
   \EO(G) = 2^{\abs{E(G)}}\pi^{-n}
      \int_{(\Rmodpi)^n} 
      \prod_{jk\in G}\cos(\theta_j-\theta_k)\,d\thetavec,
\end{equation}
where $\Rmodpi$ denotes the real numbers modulo~$\pi$.

The value of~\eqref{EOinitial2} comes mostly from a region $\varOmega_0 \subseteq (\Rmodpi)^n$ (defined in~\eqref{def:omega0}) where the
$\theta_j$ are approximately equal.
After removing the one-dimensional degeneracy caused by the invariance of the integrand under uniform translation of
each $\theta_j$, we obtain an integral of the form
\begin{equation}\label{eqxAxf}
    \int_{\varOmega} e^{-\xvec\trans\!A\xvec+f(\xvec)}\,d\xvec,
\end{equation}
where $A$ is a positive-definite matrix (defined in~\eqref{afdefs})
and $\varOmega\subseteq(\Rmodpi)^n$ is a small neighbourhood of the origin. 

Integrals of the form \eqref{eqxAxf} are common in the application of the saddle-point method for
estimating multidimensional integrals.
See Isaev and McKay~\cite{isaev2018complex} for citation of a large
number of combinatorial examples.
In \cite{isaev2018complex}, it was shown that
\[
         \int_{\varOmega} e^{-\xvec\trans\! A\xvec + f(\xvec)} d\xvec \approx \pi^{n/2} |A|^{-1/2} \Exp{ \E{ f(\X) }  + \dfrac12 
         \E{ ( f(\X) - \E{ f(\X) } )^2 } },
\]
for suitable $\varOmega$ and $f:\Reals^n\to\Complexes$,
where $\X$ is the random vector with density proportional to
$e^{-\xvec\trans\!A\xvec}$, truncated to~$\varOmega$.
In Theorem~\ref{gauss2pt}, we show that if $f(\xvec)$ 
is real and satisfies suitable
high-dimensional Lipshitz conditions in~$\varOmega$,  
integral \eqref{eqxAxf} can be approximated
to any relative precision $O(n^{-c})$ in terms of the cumulants of~$f(\X)$.

Apart from the non-symmetry of the domain, the random vector~$\X$ 
becomes a vector of independent truncated Gaussian random variables
when rotated.
This gives the more general problem of
finding the Laplace transform $\E{e^{g(\Y)}}$ of a function
of a random vector~$\Y$ with independent components.
In Theorem~\ref{cumuThm}, we provide a general expansion in terms
of cumulants when $g$ is a real function that satisfies Lipshitz conditions in the domain of~$\Y$. This relies on a new bound on the tail of the cumulant expansion, which 
is proved in Section~\ref{S:cumu-tm}.


The paper is structured as follows. 
In Section \ref{S:computations},  we review basic facts about cumulants in general and cumulants of Gaussian random variables in particular. 
We also include the tools for computing the cumulants in our formula \eqref{eo}. We also provide bounds on their rate of decrease.

In  Section \ref{S:tools}, we prove Theorem~\ref{gauss2pt}, 
and in Section~\ref{S:proof-main} we use it to prove
Theorem \ref{tm-eo} and Corollary \ref{cor-eo}.
Analysis similar to \cite{isaev2020asymptotic} is required for estimating parts of \eqref{EOinitial2} away from the region of concentration of the integral.
 
In Section \ref{S:RT}, we apply Theorem~\ref{tm-eo} to find  an asymptotic expansion for the number of regular tournaments up to the multiplicative error term of order $O(n^{-12})$.
We also obtain similar expansions for the number
of Eulerian digraphs and the number of Eulerian oriented graphs.

For $p\ge 1$, define the vector $p$-norm as
    $\norm{(x_1,\ldots,x_n)\trans}_p := \bigl(\,\sum_{i\in[n]} \abs{x_i}^p\bigr)^{\!1/p} $.
We will also use the limit as $p\to\infty$: $\norm{\xvec}_\infty:=\max_{i\in[n]} \abs{x_i}$. The same notations
are used for the induced matrix norm: $\norm{A}_p=\max_{\norm{\xvec}_p=1}\norm{A\xvec}_p$.
It is known that $\norm{A}_1$ is the maximum column sum of absolute values, and $\norm{A}_\infty$ is the maximum row sum of absolute values.

\section{Cumulants and cumulant series}\label{S:computations}

Recall that $[n]$ denotes the integer set $\{ 1, 2, \ldots, n \}$.
The combinatorial definition of joint cumulant
of random variables $\{ X_i \}_{i \in [n]}$, equivalent to \eqref{cumu},  is
\begin{equation}
\Ka{ X_1, \ldots, X_n }
= \sum_{\tau \in \mathcal P_n } (-1)^{|\tau| - 1} ( |\tau| - 1 )! \prod_{B \in \tau} 
\E{ \,\prod_{i \in B} X_i },
\label{defcumu}
\end{equation}
where $\mathcal P_n$ denotes the set of unordered partitions $\tau$ of~$[n]$ (with non-empty blocks)  and $|\tau|$ denotes  the number of blocks in the partition $\tau$. 
Recall the useful multi-linearity property of cumulants,
that is, for $r \ge 1$,
\[
\kappa\bRgg{\, \sum_{i_{1} \in [n_{1}]} X_{i_{1}}, 
\sum_{i_{2} \in [n_{2}]} Y_{i_{2}},
\ldots, \sum_{i_{r} \in [n_{r}]} Z_{i_{r}} }
= \sum_{i_{1} \in [n_{1}]} \sum_{i_{2} \in [n_{2}]} 
\cdots \sum_{i_{r} \in [n_{r}]}
\kappa\(  X_{i_{1}}, Y_{i_{2}}, \ldots, Z_{i_{r}} \),
\]
and in particular,
\[
\ka_{r}\bRgg{ \,\sum_{i \in [n]} X_{i} } 
= 
\sum_{ (i_{1}, \ldots i_{r}) \in [n]^{r}  }
\ka( X_{i_{1}}, \ldots, X_{i_{r}} ).
\]

The joint cumulant
can be regarded as a measure of the mutual dependences of the variables. 
An important property of the joint cumulant $\ka( X_1, \ldots, X_n )$ is that 
if $[n]$ can be partitioned into
two subsets $S_1$ and $S_2$ such that 
{the} variables $\{ X_i \}_{i \in S_1}$ are independent of the variables $\{ X_j \}_{j \in S_2}$, then
$\ka( X_1, \ldots, X_n ) = 0$
(see, for example, \cite{MR725217}).

\subsection{Cumulant computations for Gaussian random variables}
 
If $S$ is a set of even size, a \textit{pairing} of $S$ is a
partition of $S$ into $\card S / 2$ disjoint unordered pairs.  
Recall the following result of Isserlis~\cite{isserlis1918formula}.
\begin{lemma}\label{isserlis}
Let  $\Z=(Z_1,\ldots,Z_N)$ be a Gaussian random vector with $\E \Z = \boldsymbol{0}$ and the covariance matrix $ (\sigma_{u,v})_{u,v\in [N]}$. Consider a product $W=Z_{v_1}Z_{v_2}\cdots Z_{v_k}$, where the
subscripts do not need to be distinct.  If $k$ is odd, then
$\E W=0$.  If $k$ is even, then
\[
\E W = \sum_{\pi}  
\prod_{ij \in \pi} \sigma_{v_i, v_j},
\]
where the sum is over all pairings 
$\pi$
of $[k]$.
\end{lemma}
For example, we have
\[
\E{ Z^{2}_{1} } = \sigma_{1,1},
\quad
\E{ Z^{2}_{1} Z^{2}_{2} }  = \sigma_{11}\sigma_{2,2} + 2\sigma_{1,2}^{2},
\quad
\E{ Z_{1}^{3} Z_{2}^{3} } = 9\sigma_{1,1}\sigma_{2,2}\sigma_{1,2} + 6\sigma_{1,2}^3.
\]

In quantum field theory, 
Lemma~\ref{isserlis} is known as
Wick's formula after a later discoverer.

\begin{thm}\label{cums}
Assume the conditions of Lemma~\ref{isserlis} for an even $k$ and
let $\{P_1,\ldots,P_r\}$ be a partition of set $[k]$.
If $\pi$ is a pairing of $[k]$, define the graph $H_\pi$
as follows: $V(H_\pi)=[r]$, and for $\ell\ne m$, 
$\ell m\in E(H_\pi)$ iff $\pi$ has a pair $ij$ such that
$i\in P_\ell$ and $j\in P_m$.
Let $\varPi= \varPi(P_1,\ldots,P_m)$ be the set of all pairings $\pi$ of $[k]$ such that $H_\pi$ is connected.
Then
\[
\Ka{\,\prod_{i\in P_1} Z_{v_i},\ldots,\prod_{i\in P_r} Z_{v_i}}
=\sum_{\pi \in \varPi}  
\prod_{ij \in \pi} \sigma_{v_i v_j},
\]
where $v_1,\ldots,v_k \in  [N]$ may not necessarily be distinct.
\end{thm}
\begin{proof}
For $B\subseteq\{1,\ldots,r\}$, define $P(B)=\bigcup_{j\in B} P_j$.
From the definition of joint cumulant \eqref{defcumu}, we have that
\begin{equation}\label{cumsum}
\Ka{\,\prod_{i\in P_1} Z_{v_i},\ldots,\prod_{i\in P_r} Z_{v_i}}
=
\sum_{B_1\cup\cdots\cup B_t=\{1,\ldots,r\}}
(-1)^{t-1} (t-1)! \, \prod_{u=1}^t \,
\E{ \prod_{i\in P(B_u)} Z_{v_i} },
\end{equation}
where the sum is over all partitions of $\{1,\ldots,r\}$.
By Theorem~\ref{isserlis}, we have
\[
\prod_{u=1}^t \,\E{\, \!\prod_{i\in P(B_u)} Z_{v_i} }
=
\sum_{\pi \in\varPi_B}
\prod_{ij \in \pi} \sigma_{v_i,v_j},
\]
where $\varPi_B$ is the set of pairings 
such that no pair spans two of the sets $P(B_u)$.

Now consider any particular product
 $\prod_{ij \in \pi} \sigma_{v_i v_j}$,
 corresponding to pairing~$\pi$.
 The total weight with which this occurs in~\eqref{cumsum} is
 \begin{equation}\label{stirsum}
     \sum_{k=1}^{m} (-1)^{k-1} (k-1)! \,\stirlingii mk,
 \end{equation}
 where $m$ is the number of components of $H_\pi$
 and $\stirlingii mk$ denotes the Stirling number of the second kind (number of partitions of an $m$-set into $k$ parts).
 A standard identity for the Stirling numbers is that~\eqref{stirsum}
 equals~1 for $m=1$ and~0 for $m\geq 2$, which completes the proof.
\end{proof}

\begin{cor}\label{Cor:GeneralBound1}
   Assume the conditions of Lemma~\ref{isserlis}
    and that the covariance matrix $ \varSigma= (\sigma_{uv})_{u,v\in [N]}$ satisfies 
    $|\sigma_{u,v}|\leq 1$ for all $u,v\in [N]$.  Then, 
    for any $k_1,\ldots,k_r\in \Naturals$
    with even sum, 
    \[
       \sum_{v_1,\ldots,v_r \in [N]} |\kappa(Z_{v_1}^{k_1}, \ldots,Z_{v_r}^{k_r})|
       \leq |\varPi|\cdot \norm{\varSigma}_\infty^{r-1} N, 
    \]
    where $\varPi =\varPi(P_1,\ldots,P_r)$ is defined as in Theorem  \ref{cums}
    for some fixed  partition of 
    $[k]$ with 
    $|P_i|=k_i$ for all $i\in [r]$
    and $k = k_1+\cdots+k_r$.
    In particular, 
    $
    |\varPi| \leq (k-1)!! = (k-1)(k-3) \cdots3\cdot 1.
    $
\end{cor}
\begin{proof}
Define $\tau: [k]\rightarrow [r]$ as follows:
 for any $i \in [k]$, if  $ i\in P_j$, then $\tau(i):=j$.
 Applying Theorem \ref{cums}, we find that 
 \[
    \sum_{v_1,\ldots,v_r \in [N]} |\kappa(Z_{v_1}^{k_1}, \ldots,Z_{v_r}^{k_r})|
    \leq \sum_{\pi \in \varPi}
    \sum_{v_1,\ldots,v_r \in [N]} \prod_{ij \in \pi} |\sigma_{v_{\tau(i)},v_{\tau(j)}}|.
 \]
It is sufficient to show that, for any $\pi \in \varPi$, 
 \[
 \sum_{v_1,\ldots,v_r \in [N]} \prod_{ij \in \pi} |\sigma_{v_{\tau(i)},v_{\tau(j)}}| \leq 
   \norm{\varSigma}_\infty^{r-1} N, 
 \]
 We prove it by induction on $r$. If $r=1$, then by the assumption that $|\sigma_{u,v}|\leq 1$ for all $u,v\in [N]$, 
 we have
  \[
   \sum_{v_1,\ldots,v_r \in [N]} \prod_{ij \in \pi} |\sigma_{v_{\tau(i)},v_{\tau(j)}}|
   = 
   \sum_{v_1 \in [N]} 
   |\sigma_{v_{1},v_{1}}|^{\abs{\pi}}
   \leq N.
  \]
For the induction step, assume $r\geq 2$ and  observe that 
 since $H_{\pi}$ is connected there is a vertex in $V(H_{\pi})$ such that the graph remains connected if we remove it from $H_{\pi}$. Without loss of generality, we can assume that this vertex is $r$ adjacent to vertex $x$. Let $\pi'$ be the pairing obtained from $\pi$ by removing all pairs containing a point from the part $P_r$ corresponding to vertex $r$ of $H_\pi$. Estimating other covariances by 1, we find that
\[
        \sum_{v_r\in [N]} 
  \prod_{ij\in \pi\setminus \pi'} |\sigma_{v_{\tau(i)},v_{\tau(j)}}| \leq \sum_{v_r\in [N]} 
  |\sigma_{v_{x},v_{r}}| \leq \norm{\varSigma}_\infty.
\]
Using the induction hypothesis, we have 
 \[
  \sum_{v_1,\ldots,v_r \in [N]} \prod_{ij \in \pi} |\sigma_{v_{\tau(i)},v_{\tau(j)}}| \leq
    \sum_{v_1,\ldots,v_{r-1}\in [N]}  \norm{\varSigma}_\infty\prod_{ij\in \pi'} |\sigma_{v_{\tau(i)},v_{\tau(j)}}|
    \leq  \norm{\varSigma}_\infty^{r-1} N. 
 \]
 This proves the induction step and completes the proof of the corollary.
\end{proof}

\subsection{Bounds on cumulants in Theorem \ref{tm-eo}} 
Recall that, for   a connected graph $G$ on vertex set $[n]$,   $\X_G$ denotes  $n$-dimensional singular Gaussian  random vector on
 the subspace $\mathcal{V}_n$ defined in~\eqref{eq:Vn}
with density proportional to 
$\Exp{-\frac12\xvec\trans\! L\xvec}$, where $L=L(G)$ is the Laplacian matrix of graph $G$ defined by \eqref{Ldef}.
Recall from \eqref{eq:eoTaylor} that 
 \[
f_{K}(\x)
= \sum_{t=2}^K c_{2t} \sum_{jk \in G} (x_j-x_k)^{2t}.
\]
From \cite[eqn.~23.1.14]{Abramowitz},
 we know that 
$\abs{B_{2t}}< 4 (2t)!/(2\pi)^{2t}$,
and therefore, by \eqref{eq:cDef},
\begin{equation}\label{c-bound}
  | c_{2t} |
  <
  \frac{2^{2t+1} (4^{t}-1)}
  {t (2\pi)^{2t}}
  <
  \Dfrac{2}{t}
  \left( \Dfrac{2}{\pi} \right)^{2t}.
\end{equation}

In the next theorem, we prove bounds on the cumulants of $f_K(\X_G)$ in terms of the infinity norm of a 
non-singular Gaussian random variable $\X$ that projects into $\X_G$.  
Let $J_n$ denote the matrix with every entry one.
\begin{lemma}\label{L:projections}
Let $L= L(G)$ be the Laplacian matrix of a connected graph $G$.
 Let  $\varSigma_w$ be defined by $\varSigma_w^{-1} = L+wJ_n$ for some $w>0$. If  $\X$  is a   Gaussian  random vector  
with density  
$(2\pi)^{-n/2} |\varSigma_w|^{-1/2}\exp\bigl(-\frac{1}{2}\xvec\trans \varSigma_w^{-1}\xvec\bigr)$ then 
$
g(\X_G) \stackrel{d}{=} g(\X) 
$
for any function $g$ such that $g(\xvec+ \theta \boldsymbol{1}) = g(\xvec)$ for all $\xvec \in \Reals^n$
and $\theta\in\Reals$.
\end{lemma}
\begin{proof}
Since $G$ is connected, $L+wJ_n$ is 
positive-definite for all $w>0$.
Note that $I-J_n/n$ is the projector operator to the space $\mathcal V_n$.
Applying this linear transformation to the Gaussian vector $\X$, we find that
$(I-J_n/n)\X$ has the same distribution as $\X_G$.
Therefore, $g(\X_G) \stackrel{d}{=} g((I-J_n/n)\X) = g(\X)$ as claimed.
\end{proof}

\begin{thm}\label{T:cumulant_bounds}
Let $G$ be a connected graph with
maximum degree $\dmax$, and minimum degree~$\dmin$.  
Assume that, for some $w \in (0,1]$,
\[
    \varSigma_w^{-1} = L  + wJ_n, \qquad \norm{\varSigma_w}_\infty \leq 1/2,
\]
where $L= L(G)$ is the Laplacian  matrix of $G$. Then the following bounds hold.
\begin{itemize}
\item[(a)]
 For  $K \leq \dmin/2$, we have
\[
    \E {f_K(\X_G)} 
= -\Bigl(\dfrac14 + O(\dmin^{-1})\Bigr) \sum_{ jk \in G } 
\bR{ \Dfrac{1}{d_{j}}
+ \Dfrac{1}{d_{k}} }^2
+ \O{ \Dfrac{n}{\dmin} \infnorm{\varSigma_w} + \Dfrac{n\dmax}{\dmin^2} \infnorm{\varSigma_w}^2},
\]
where $d_1,\ldots,d_n$ are the degrees of $G$.
\item[(b)]
For $r,K$ such that $r(K+2) \leq \dmin/4$, we have   
\[
   |\kappa_{r}(f_K(\X_G))| \leq  
   \Dfrac{n}{2\dmin} \left(\Dfrac{5\dmax}{\dmin}\right)^{r} \norm{\varSigma_w}_\infty^{r-1}    (4r-1)!!\,.
\]
\item[(c)]
For $r, K, K'$   such that $K'>K$
and $(r+1)(K'+2) \leq \dmin/4$, we have
\[
   |\kappa_{r}(f_K(\X_G))- \kappa_{r}(f_{K'}(\X_G))| \leq 
    \Dfrac{r n}{2\dmin} \left(\Dfrac{5\dmax}{\dmin}\right)^{r} \left(\Dfrac{2}
    {\dmin}\right)^{K-1}\norm{\varSigma_w}_\infty^{r-1}    (2K+4r-3)!!\,.
\]
\end{itemize}
\end{thm}
\begin{proof}
Recalling  the definition of $f_K$ and 
 applying Lemma \ref{L:projections}, we find that 
 \[
    f_K(\X_G)
    \stackrel{d}{=} f_K(\X),
 \]
where $\X$ is a   Gaussian  random vector  
with density  
$(2\pi)^{-n/2} |\varSigma|^{-1/2}\Exp{-\frac12\xvec\trans \varSigma_w^{-1\,}\xvec}$.
Thus, it is sufficient to prove the bounds for  $\E{f_K(\X)}$ and $\kappa_r(f_K(\X))$.

First, we  estimate the entries of  the matrix 
  $\varSigma_w = (\sigma_{jk})$.
Define the diagonal matrix  $D$  by  $D_{jj} = d_j$, where $d_j$ is the degree of vertex $j$ in $G$.
Observing also that all entries of $\varSigma_w^{-1} - D= L+wJ - D$ are in $[-1,1]$, 
we obtain  
\begin{align*}
\norm{\varSigma_w - D^{-1}}_{\max}
&=\norm{\varSigma_w(D- \varSigma_w^{-1}) D^{-1}}_{\max} \\
&\le 
\infnorm{\varSigma_w}
\norm{(D-\varSigma_w^{-1})D^{-1}}_{\max}
\leq \infnorm{\varSigma_w} \max_{j \in [n]} \Dfrac{1}{d_j} \leq \Dfrac{1}{\dmin} \infnorm{\varSigma_w},
\end{align*}
where $\norm{\cdot}_{\max}$ denotes the maximum absolute value of entries.  This gives
\begin{equation}\label{eq:sigma_bounds}
\left|\sigma_{jj} - \dfrac{1}{d_j}\right|
            \leq \Dfrac{1}{\dmin}\,\infnorm{\varSigma_w}, 
            \qquad |\sigma_{jk}|\leq \Dfrac{1}{\dmin} \,\infnorm{\varSigma_w}.
\end{equation}

Let 
$X_{jk} := X_j - X_k$
for distinct $j, k \in [n]$.
Let 
\begin{equation}\label{sigma-jkst}
\sigma_{jk, st} = \Cov( X_{jk}, X_{st} )
= \sigma_{js}- \sigma_{jt} -\sigma_{ks}+\sigma_{kt}.
\end{equation}
Using  \eqref{eq:sigma_bounds}, for all distinct $j,k\in [n]$,
\[
    \sigma_{jk,jk} = \Dfrac{1}{d_j} + \Dfrac{1}{d_k} 
    \pm \Dfrac{4}{\dmin}\infnorm{\varSigma_w}, 
\]
where we mean that the last term lies in
the interval
$\bigl[-\dfrac{4}{\dmin}\infnorm{\varSigma_w},\dfrac{4}{\dmin}\infnorm{\varSigma_w}\bigr]$.
Then, by the assumption that   $\norm{\varSigma_w}_\infty\leq 1/2$, we obtain
\begin{equation}\label{eq:sigma-sigma}
|\sigma_{jk, st}|
\leq 
\begin{cases}
   \,\dfrac{4}{\dmin}, & 
    \text{if } 
    \abs{\{ j, k \} \cap  \{ s, t \}}
    \geq  1; \\[1ex]
    \,\dfrac{4}{\dmin}\infnorm{\varSigma_w}
    , & \text{if } 
    \abs{\{ j, k \} \cap  \{ s, t \}}
    = 0.
\end{cases}
 \end{equation}
 From Lemma \ref{isserlis}, we know that
 \[
 \E{ X_{jk}^{2t} } 
 = (2t-1)!!\, \sigma_{jk}^{t}.
 \]
Recalling 
the definition of $f_K(\xvec)$ from \eqref{eq:eoTaylor},
and using \eqref{c-bound}
and \eqref{eq:sigma_bounds}, 
we obtain that 
\begin{align*}
\E {f_K(\X)} 
&= \sum_{t=2}^K c_{2t}
\sum_{jk \in G}
\E{
X_{jk}^{2t} } 
= -\Dfrac{1}{12}
\sum_{jk \in G}
\E{ X_{jk}^{4} }
+ \sum_{t=3}^K c_{2t}
\sum_{jk \in G}
\E{
 X_{jk}^{2t} }
 \\
 &=  -\Dfrac14 
\sum_{ jk \in G } 
\bR{ \Dfrac{1}{d_{j}}
+ \Dfrac{1}{d_{k}} \pm \Dfrac{4}{\dmin} \norm{\varSigma_w}_\infty}^2 
\pm
n\dmax \sum_{t=3}^K \bR{ \Dfrac{2}{\pi} }^{2t} (2t-1)!! \bR{ \Dfrac{4}{\dmin} \norm{\varSigma_w}_\infty }^{t}.
\end{align*}
Since $x(2t-x)\leq t^2$, we can  bound $(2t-1)!! \leq t^{t}$. Recalling the assumption that $K \leq \dmin/2$, we estimate 
 \[
    n\dmax\sum_{t=3}^K 
\bR{ \Dfrac{2}{\pi}
 }^{2t}
(2t-1)!!\,
\bR{ \Dfrac{4}{\dmin}
\norm{\varSigma_w}_\infty}^{t}  = \O{\Dfrac{n\dmax}{\dmin^3}
\norm{\varSigma_w}_\infty^3}. \]
We also have 
\[
    \sum_{ jk \in G } 
\bR{ \Dfrac{1}{d_{j}}
+ \Dfrac{1}{d_{k}} \pm \Dfrac{4}{\dmin} \norm{\varSigma_w}_\infty}^2  = 
\sum_{ jk \in G } 
\bR{ \Dfrac{1}{d_{j}}
+ \Dfrac{1}{d_{k}} }^2  \pm \Dfrac{16 n }{\dmin}\norm{\varSigma_w}_\infty \pm  \Dfrac{16 n\dmax}{\dmin^2}\norm{\varSigma_w}_\infty^2.
\]
This proves (a).

We proceed to part (b).
Using \eqref{sigma-jkst},
for any distinct $j, k \in [n]$, we have that
\begin{equation}\label{sum:sigma-sigma}
\sum_{st \in G}
\Abs{\sigma_{jk, st}}
\le
2\sum_{st \in G}  
\Abs{ \sigma_{j s} }
+2\sum_{st \in G}  
\Abs{ \sigma_{k s} }
\le
4\dmax
\infnorm{\varSigma_w}.
\end{equation}
By linearity,
\begin{equation}
\begin{aligned}
\label{cum-expan}
\kappa_r \bR{ f_{K}(\X) }
&= 
\kappa_r \bRgg{\,
\sum_{t=2}^K c_{2t} \sum_{jk \in G} (X_j-X_k)^{2t} } \\
&= 
\sum_{
 t_1 = 2}^K \cdots\sum_{
 t_r = 2}^K\
\sum_{ e_1,\ldots,e_r \in  G}
\bRgg{ \,\prod_{s \in [r]} c_{2t_s} }
\kappa \bR{
X_{e_1}^{2t_1},
\ldots,
X_{e_r}^{2t_r} }.
\end{aligned}
\end{equation}
Applying Corollary \ref{Cor:GeneralBound1} for 
$(Z_v)_{v\in [N]}:=\dfrac{\dmin^{1/2}}{2}(X_{jk})_{jk \in G}$ and
$N=\abs{E(G)}\le n\dmax/2$.
From \eqref{eq:sigma-sigma}, we have that the covariances
of $(Z_v)_{v\in [N]}$ are bounded by $1$ as required.  
Using \eqref{sum:sigma-sigma} to bound the infinity norm of the  covariance matrix, we find that 
    \begin{equation}\label{sum-edges-cum}
\sum_{ e_1,\ldots,e_r \in  G}
|\kappa \bR{
X_{e_1}^{2t_1},
\ldots,
X_{e_r}^{2t_r} }|
\leq  \Dfrac{n\dmax}{2} \left(\Dfrac{4}{\dmin}\right)^{t} 
(\dmin \dmax\norm{\varSigma_w}_\infty)^{r-1} (2t-1)!!, 
\end{equation}
where $t = t_1+\cdots+t_r$. 
Since $t_i \in \{2,\ldots,K\}$,  
using $(2r+t-x)(2r+t+x)
\le (2r+t)^2$,
we have
\[
    (2t-1)!! \leq (2r+t)^{t -2r}  (4r-1)!!\leq  
     ((K+2)r)^{t -2r} (4r-1)!!.
\]
Recalling our assumption $r(K+2) \leq \dmin/4$,
and
substituting the above bounds and \eqref{c-bound} into \eqref{cum-expan}, we obtain  
\begin{align*}
    |\kappa_r(f_K(\X))|  &\leq 
    \Dfrac{n\dmax}{2}(\dmin \dmax \norm{\varSigma_w}_\infty)^{r-1}
    (4r-1)!!
    \sum_{
 t_1 = 2}^K \cdots\sum_{
 t_r = 2}^K  \left(\Dfrac{2}{\pi}\right)^{2t} \left(\Dfrac{4}{\dmin}\right)^{t} 
   ((K+2)r)^{t -2r}  
   \\ 
   &\leq \Dfrac{n\dmax}{2}(\dmin \dmax \norm{\varSigma_w}_\infty)^{r-1}
    \left(\Dfrac{16}{\pi^2 \dmin}\right)^{2r}  (4r-1)!!
    \left(\sum_{i =0}^{K-2} \left(\Dfrac{16 r(K+2)}{\pi^2 \dmin}\right)^{\!i\,} \right)^r
    \\& \leq  \Dfrac{n}{2\dmin}\left(\Dfrac{\dmax}{\dmin}\right)^{\!r} \norm{\varSigma_w}_\infty^{r-1}    (4r-1)!! 
     \left(\Dfrac{256}{\pi^4}\sum_{i =0}^{\infty}   \Bigl(\Dfrac{4}{\pi^2}\Bigr)^{i} \right)^{\!r}.
\end{align*}
Computing $\dfrac{256}{\pi^4}\sum_{i=0}^{\infty}   \Bigl(\dfrac{4}{\pi^2}\Bigr)^{i} \approx 4.42 < 5$, part (b) follows.

For part (c),  using \eqref{c-bound}, \eqref{cum-expan} and \eqref{sum-edges-cum}, we estimate 
 \begin{equation}
\begin{aligned}
\label{cum-expan2}
\bigl|\kappa_r \bR{ f_{K}(\X) } &- \kappa_r \bR{ f_{K'}(\X) }\bigr|
\\
&\leq  
r\sum_{
  t_1 = K+1}^{K'} \sum_{
 t_2 = 2}^{K'}\ \cdots\sum_{
 t_r = 2}^{K'}\
\sum_{ e_1,\ldots,e_r \in  G}
\bRgg{ \,\prod_{s \in [r]} c_{2t_s} }
\bigl|\kappa \bR{
X_{e_1}^{2t_1},
\ldots,
X_{e_r}^{2t_r} }\bigr|
\\
&\leq  r \sum_{
  t_1 = K+1}^{K'} \sum_{
 t_2 = 2}^{K'}\ \cdots\sum_{
 t_r = 2}^{K'}
 \Dfrac{nd}{2} \left(\Dfrac{16}{\dmin \pi^2}\right)^{t} 
(\dmin \dmax \norm{\varSigma_w}_\infty)^{r-1} (2t-1)!!,
\end{aligned}
\end{equation}
where $t = t_1 + \cdots + t_r$ as before.
Similarly, we obtain
\[
        (2t-1)!! \leq (t + K+2r-1)^{t - K -2r+1} (2K+4r-3)!!
        \leq ((r+1) (K'+2))^{t - K -2r+1}(2K+4r-3)!! .
\]
Then, using $ (r+1) (K'+2) \leq \dmin/4$, we can bound
\begin{align*}
 \sum_{
  t_1 = K+1}^{K'} \sum_{
 t_2 = 2}^{K'}\ \cdots&\sum_{
 t_r = 2}^{K'}
   \left(\Dfrac{16}{\dmin \pi^2}\right)^{t} 
  (2t-1)!!
  \\
  &\leq \left(\Dfrac{16}{\dmin \pi^2}\right)^{K+2r-1} (2K+4r-3)!!
  \left(\sum_{i =0}^{K'}
   \left(\Dfrac{16 (r+1) (K'+2)}{\dmin \pi^2}\right)^{\!i\,}
  \right)^{\!r}
  \\
  &\leq  5^r \left(\Dfrac{16}{ \pi^2 \dmin }\right)^{K-1} \dmin^{-2r}
  (2K+4r-3)!!
\end{align*}
Substituting this bound into \eqref{cum-expan2} completes the proof.
\end{proof}



\section{Cumulant expansion of multidimensional integrals}
\label{S:tools}

In this section, we introduce in Theorem~\ref{cumuThm} a new explicit tail bound
on the cumulant expansion of the Laplace transform
of a real function of independent random variables.
Then, in Theorem~\ref{gauss2pt} we apply Theorem~\ref{cumuThm}
to integrals of the form
\[
   \int_{\varOmega} e^{-\xvec\trans\! A\xvec + f(\xvec)}\,d\xvec
\]
for suitable domains $\varOmega$ and functions
$f:\Reals^n\to\Reals$, where $A$ is a positive-definite
$n\times n$ matrix.
Such integrals appear in applications of the saddle-point
method to asymptotics of integrals.

\subsection{Laplace transform of real functions of independent random variables}

Let $\X=(X_1,\ldots, X_n)$ be a random vector with independent components
taking values in 
a product space $\S := S_1\times\cdots\times S_n$.
For $\yvec \in \S$
and $f : \S \to \Reals$,  let 
$R^j_{\yvec} $ denote the operator that
replaces $f(x_1,\ldots x_n)$
by $f(x_1,\ldots, x_{j-1}, y_j, x_{j+1},\ldots x_n)$:
\[
 R^j_{\yvec}  [f] (\xvec) =  f(x_1,\ldots, x_{j-1}, y_j, x_{j+1},\ldots x_n). 
\]
For  $V = \{v_1,\ldots, v_k\}\subseteq [n]$, define 
\begin{equation}\label{def:R-delta}
R^{V}_{\yvec}= R^{v_1}_{\yvec}\cdots R_{\yvec}^{v_k}
\qquad 
\mbox{ and }
\qquad
\partial^V_{\yvec}:= \partial_{\yvec}^{v_1}\cdots \partial_{\yvec}^{v_k},
\end{equation}
where  $\partial_{\yvec}^j:= I - R^j_{\yvec}$ and  $I$ is the identity operator.  
Let
\begin{equation}\label{def:Delta}
	 \Delta_{V}(f)
=	 \Delta_{V}(f, \S)  :=  \sup_{ \xvec,\y\in \S} \left|\partial^V_{\yvec}[f] (\xvec)\right|.
\end{equation}
Note that these operators commute, so the order of
elements in $V$ doesn't matter.

\begin{thm}\label{cumuThm}
	Let $m > 0$ be an integer and $\alpha \geq 0$. 
	Suppose $f$ is a bounded real function on $\S$ 
	  such that 
\[
	 \max_{\substack{v \in [m]\\ j \in [n]}}\,\sum_{V\in \binom{[n]}{v} \st j \in V } \!\Delta_V(f) \leq \alpha.
\]
	Then, we have
   \[
  	 \E { e^{f(\X)} }  =  
   (1+\delta)^n  \exp\left(\sum_{r=1}^m \frac{\kappa_r (f(\X))}{r!} \right),
  \]
   where $\delta> -1$ satisfies $\abs\delta \leq  e^{(100 \alpha)^{m+1}}-1$. 
   Furthermore,
   for any $r\in [m]$,
   \[
            | \kappa_r (f(\X))| 
            \leq  n \Dfrac{(r-1)!}{50r}(80
             \alpha)^{r}.  
   \]
\end{thm}

We will prove Theorem \ref{cumuThm} in Section \ref{S:cumu-tm}.
The case $m=2$ is similar to the second-order approximation of the Laplace transform considered by Catoni in \cite{catoni2003laplace}.
A significant improvement 
in comparison with \cite[Theorem~1.1.]{catoni2003laplace} 
(apart from having better precision) is that we only require averages of $\Delta_V(f)$
to be small while some of them can be quite big.
The case of $m=2$ was extended to complex functions in~\cite{isaev2018complex}, but note that Theorem~\ref{cumuThm} is only for real functions.

\subsection{Moments and cumulants of truncated Gaussians}

In this section we prove a lemma
for  approximating
cumulants of truncated Gaussian 
by their values for the unrestricted Gaussian. 
We will need the following two lemmas.
\begin{lemma}
\label{prodDif}
If we have that $a_1, a_2, \ldots, a_n, b_1, b_2, \ldots, b_n \in [-1, 1]$,
then
\begin{equation}\label{diffprod}
\biggl| \prod_{i \in [n]} a_i
- \prod_{i \in [n]} b_i \biggr|
\le \sum_{i \in [n]} \Abs{ a_i - b_i }.
\end{equation}
\end{lemma}
\begin{proof}
We prove it by induction on $n$.
The base case when $n=1$ holds trivially.
Suppose \eqref{diffprod} holds for some $n \ge 2$.
Let $A_n$ and $B_n$ denote 
$\prod_{i \in [n]} a_i$ and $\prod_{i \in [n]} b_i$ respectively.
Then 
\begin{align*}
\left| A_{n+1}
- B_{n+1} \right|
= \left|  a_{n+1}  A_n
- b_{n+1} B_n \right|
&= \left|  (a_{n+1} - b_{n+1}) A_n
- b_{n+1} (A_n - B_n) \right| \\
&\le \left|  a_{n+1} - b_{n+1} \right|
+  \left| A_n - B_n \right| 
\le \sum_{i \in [n+1]} \Abs{ a_i - b_i }.
\end{align*}
This completes the proof.
\end{proof}

\begin{lemma}\label{ordpart}
Let $N_{s,k}$ be the number of ordered partitions of
$[s]$ with $k$ parts. Then
\[
  N_{s,k} \le 2^{k-s}\binom{s-1}{k-1} s!.
\]
\end{lemma}
\begin{proof}
We can bound $N_{s,k}$ as follows.  Take a permutation
of $s$ and cut it into $k$ non-empty parts (e.g., for $s=5,k=3$
we could have $[3|4,1|2,5]$).
The cutting can be done in $\binom{s-1}{k-1}$ ways, but there
is an overcount due to the ordering
within the parts, so we can weight each cutting by $(b_1!\cdots b_k!)^{-1}$
where $b_1,\ldots,b_k$ are the part sizes.
Using $b_j!\ge 2^{b_j-1}$, this gives the desired bound.
\end{proof}

Now, we are ready to prove the lemma  relating truncated and non-truncated cumulants.

\begin{lemma}\label{truncated}
Let $A$ be an $n\times n$ symmetric positive-definite real matrix.
Let $\X:\Reals^n\to\Reals$ be a random vector with density
\[
    \pi^{-n/2} \abs{A}^{1/2} e^{-\xvec\trans\! A\xvec}.
\]
Let $\varOmega\subseteq\widehat\varOmega\subseteq\Reals^n$ be measurable
sets and define $p=\pr(\X\notin \varOmega)\le\frac34$ and $\hat p=\pr( \X\notin\widehat\varOmega )$.
Let $f:\Reals^n\to\Reals$ be measurable such that
for some $b, \hat b \ge 0$,
\[
   \abs{f(\xvec)} \le \begin{cases}
       \,Ce^{b\xvec\trans\! A\xvec/n}, & \text{for $\xvec\in\Reals^n$}; \\[1ex]
       \,\widehat Ce^{\hat b\xvec\trans\! A\xvec/n}, & \text{for $\xvec\in\widehat\varOmega$}.
   \end{cases}
\]
Then, if 
$n\ge \max\{br+b^2r^2,\hat br+\hat b^2r^2\}$,
\[
\Abs{ \kappa_r \(f(\X) \mid \X\in\varOmega\) - \kappa_r (f(\X)) } 
\le 6^{r+1} r!
\bigl(
C^r e^{br/2} \hat p^{1-br/n} + 
\widehat C^r e^{\hat br/2} (p^{1-\hat br/n}+\hat p^{1-\hat br/n})\bigr).
\]
\end{lemma}

\begin{proof}
Define
\[
   \hat f(\xvec) = \begin{cases}      
       f(\xvec), & \text{for $\xvec\in\widehat\varOmega$}; \\
         0, & \text{for $\xvec\notin\widehat\varOmega$}.
   \end{cases}
\]
Then by the triangle inequality
and the assumption
$\varOmega\subseteq\widehat\varOmega$, we have 
\begin{equation}\label{triangle}
\begin{aligned}
\Abs{ \kappa_r \(f(\X) &\mid \X\in\varOmega\) - \kappa_r (f(\X))}
\le
\Abs{ \kappa_r \(f(\X) \mid \X\in\widehat\varOmega\) - \kappa_r (f(\X))} \\
&\quad+ 
\Abs{ \kappa_r \( \hat f(\X) \mid \X\in\varOmega\) - \kappa_r ( \hat f(\X )} 
+ \Abs{ \kappa_r \( \hat f(\X) \mid \X\in\widehat\varOmega \) - \kappa_r ( \hat f(\X))}.
\end{aligned}
\end{equation}

\smallskip

The three terms on the right have the same structure, so we will continue with a generic case.
Suppose $g(\xvec)$ satisfies $\abs{g(\xvec)}\le e^{\tilde b\xvec\trans\! A\xvec/n}$ for all $\xvec$, and $\tilde p=\pr(\X\notin\medtilde\varOmega)\le\frac34$.
From \eqref{defcumu}, we have
\begin{align}
\kappa_r ( g(\X) &\mid \X\in\medtilde\varOmega) - \kappa_r ( g(\X)) \notag \\
&= \sum_{\tau \in \mathcal P_r} (-1)^{|\tau| - 1} ( |\tau| - 1 )! 
\bRg{\, \prod_{B \in \tau} 
\EE{  g(\X)^{|B|} } - \prod_{B \in \tau} 
\EE{  g(\X)^{|B|} \mid \X\in\medtilde\varOmega } }.
\label{cumu-gau-trun}
\end{align}

Let $B\in\tau\in \mathcal P_r$.
Applying \cite[Lemma 4.1]{isaev2018complex} to $g(\X)^{|B|}$,
we have
\[
\Abs{ \EE{  g(\X)^{|B|} }
-  \EE{  g(\X)^{|B|} \mid\X\in\medtilde\varOmega } }
\le 15\, e^{b|B|/2} \tilde p^{1-b|B|/n}.
\]
Furthermore, we can bound
\[
\left| \E{ g(\X)^{|B|} } \right|
\le \E{ \Abs{ g(\X)^{|B|} } }
\leq \frac{\int_{\Reals^n} e^{-  (1-\frac{b|B|}{n})\xvec\trans\! A \xvec} d\xvec}{\int_{\Reals^n} e^{-   \xvec\trans\! A \xvec} d\xvec}
 = (1 - b|B|/n)^{-n/2},
\]
and similarly, using 
$\EE{ \Z \mid \medtilde\varOmega } 
\le \EE{ \abs{\Z} } / \pr(\medtilde\varOmega)$ for any random variable~$\Z$,
we have
\[
\Abs{\EE{ g(\X)^{|B|} \mid\X\in\medtilde\varOmega } }
 \le (1-\tilde p)^{-1}(1-b|B|/n)^{-n/2}
 \le 4(1-b|B|/n)^{-n/2}.
\]
For $x\in(0,1)$, Taylor expansion shows that 
$x<-\log(1-x)<x + \dfrac{x^2}{2(1-x)}$.
Since $\abs{B}\le r$ and $n\ge br+b^2r^2$, we have
\[
e^{b|B|/2} \le (1-b|B|/n)^{-n/2} \le e^{b|B|/2+1/4}.
\]
In order to apply Lemma~\ref{prodDif} to estimate each term of \eqref{cumu-gau-trun},
we scale $\EE{  g(\X)^{|B|} }$ and $\EE{  g(\X)^{|B|} \mid \X\in\medtilde\varOmega }$ by $4(1-b|B|/n)^{-n/2}$.
For any $\tau\in \mathcal P_r$, this gives
\begin{align*}
\biggl| \prod_{B \in \tau} \EE{  g(\X)^{|B|} }
&- \prod_{B \in \tau} \EE{  g(\X)^{|B|} \mid \X\in\medtilde\varOmega } \biggr|
\\ &\le
\bR{ \prod_{B\in\tau} \bigl(4(1-b|B|/n)^{-n/2}\bigr) }
\cdot
\bR{ \sum_{B\in\tau} \Dfrac{15e^{b|B|/2}\tilde p^{1-b|B|/n}}{4(1-b|B|/n)^{-n/2}}
}
 \\ &\le 
 4^{|\tau|+1} \prod_{B\in\tau} e^{b|B|/2+1/4} \cdot
  \sum_{B\in\tau} \tilde p^{1-b|B|/n} 
 \\ &\le
   4^{\abs{\tau}+1} \abs{\tau}\,e^{br/2+\abs{\tau}/4} \tilde p^{1-br/n}.
\end{align*}
Consequently, from~\eqref{cumu-gau-trun}, we have
\begin{align*}
  \Abs{\kappa_r ( g(\X) \mid \X\in\medtilde\varOmega) - \kappa_r (g(\X))}
  &\le 4 e^{br/2} \tilde p^{1-br/n}
  \sum_{\tau\in \mathcal P_r} 4^{\abs{\tau}} e^{\abs{\tau}/4} \,\abs{\tau}! \\
  &= 4 e^{br/2} \tilde p^{1-br/n} 
   \sum_{k=1}^r 4^k e^{k/4} N_{r,k}, 
\end{align*}
where $N_{r,k}$ is the number of 
ordered partitions of $[r]$ with $k$ parts.
Applying Lemma~\ref{ordpart},
\begin{align*}
\bigl|\kappa_r (g(\X) \mid \X\in\medtilde\varOmega) - \kappa_r (g(\X))\bigr|
&\le 4 r!\,e^{br/2} \tilde p^{1-br/n} 2^{-r}
   \sum_{k=1}^r \binom{r-1}{k-1} 8^k e^{k/4} \\
&= 4 r!\,e^{br/2} \tilde p^{1-br/n} 2^{-r} (1+8e^{1/4})^{r-1} 8e^{1/4} \\
& \le 6^{r+1} r!\,e^{br/2} \tilde p^{1-br/n},
\end{align*}
as desired.
Note that $\kappa_r(c\Z)=c^r\kappa_r(\Z)$ for any random variable~$\Z$
and constant~$c$.
Applying the corresponding bounds to each of the terms in the
right side of \eqref{triangle}, we obtain
\begin{align*}
 \Abs{ \kappa_r \(f(\X) &\mid \X\in\varOmega\) - \kappa_r (f(\X))} \\
&\le 
6^{r+1} r!
\bR{
C^r e^{br/2} \hat p^{1-br/n}
+ 
\widehat C^r e^{\hat br/2} p^{1-\hat br/n}
+ 
\widehat C^r e^{\hat br/2} \hat p^{1-\hat br/n}
}.
\end{align*}
This completes the proof.
\end{proof}

\subsection{Laplace transforms of functions of truncated Gaussian random vectors}

Recall that our aim is to estimate integrals of the type  
\[
\int_\varOmega e^{-\xvec\trans\!A\xvec + f(\xvec)}\,d\xvec
 = \E{e^{f(\X)}} \int_{\varOmega} e^{-\xvec\trans\! A\xvec}\,d\xvec,
\]
where $\X$ is a truncated  Gaussian random vector with support
$\varOmega$ and density proportional to $e^{-\xvec\trans\!A\xvec}$ and $f:\Reals^n\to\Reals$.
To apply Theorem \ref{cumuThm}, we consider $g(\yvec):=f(T^{-1}\xvec)$, where $T$ is a linear transformation such that $T\trans\!A T=I$.
Then $\Y =T^{-1}\X$ is a standard
Gaussian vector with independent components, truncated to~$T^{-1}\varOmega$.

We can bound the quantities $\Delta_V(g)$ required in Theorem \ref{cumuThm} in terms of the mixed derivatives of $f$, using the following lemma.
For $\varOmega \subseteq \Reals^n$
and some analytic $f: \varOmega \to \Reals$,
for $\ell\in[n]$, define 
  \begin{equation}\label{def:delta-bar}
     \medbar\Delta_\ell(f, \varOmega) := 
        \max_{u^{}_1\in[n]} 
        \sum_{u^{}_2,\ldots,u^{}_\ell\in [n]}\,
          \sup_{\xvec\in \Omega}\;
        \biggl| \Dfrac {\partial^\ell f(\xvec)}{\prod_{r\in [\ell]}\partial x_{u^{}_r}}\biggr|.
  \end{equation}
 For $\rho \geq0$,  let
 \[
        U_n(\rho):= \{\xvec \in \Reals^n \st \norm{\xvec}_\infty \leq \rho\}. 
 \]

\begin{lemma}\label{Deltabound}
  Suppose 
   $g :  U_n(\rho) \to \Reals$ is an analytic function. 
  Let 
  $f : \varOmega \to \Reals$ be defined by
  $f(\xvec):=g(T^{-1}\xvec)$, where $\varOmega = T ( U_n(\rho))$ and 
   $T$ is a real $n\times n$ invertible matrix.
  Then for $\ell\in[n]$, we have
  \[
      \max_{j\in [n]} \sum_{V\in\binom{[n]}{\ell}:j\in V} \Delta_V(g)   
      \leq 
      \frac{\infnorm{T}^{\ell-1}\, \onenorm{T}}
 {(\ell-1)!}
 (2\rho)^{\ell}
         \medbar\Delta_\ell(f, 
         \varOmega),
  \]
  where 
  $\Delta_V(g) = \Delta_V(g, U_n(\rho))$
  is defined according to \eqref{def:Delta},
  and
     $\medbar\Delta_\ell(f, \varOmega)$
is defined by \eqref{def:delta-bar}.  
\end{lemma}

\begin{proof}
 Applying the mean value theorem, we obtain
  \[
      \Delta_V(g) \leq 
         \sup_{\yvec\in U_n(\rho)}\;\biggl| \Dfrac{\partial^{\ell} g(\yvec)}{\prod_{r\in V}\partial y_r} \biggr| 
            \, (2\rho)^\ell.
  \]
Then it suffices to show that
\[
     \widehat\Delta_\ell(g) := \max_{j\in [n]} \sum_{V\in\binom{[n]}{\ell}:j\in V}\,
       \sup_{\yvec\in U_n(\rho)}\;
        \biggl| \Dfrac{\partial^\ell g(\yvec)}{\prod_{r\in V}\partial y_r} \biggr|
  \]
satisfies the bound
  \[
     \widehat\Delta_\ell(g)
        \leq \frac{\infnorm{T}^{\ell-1}\, \onenorm{T}}
 {(\ell-1)!}
      \; \medbar\Delta_\ell(f, 
         \varOmega ).
  \]

  Let $T=(t_{jk})$.
    Since $\xvec=T\yvec$, by the chain rule,
    for any distinct $j_1,\ldots,j_\ell\in[n]$, we have that
    \[
       \Dfrac{\partial^\ell g(\yvec)}{\partial y_{j_1}\cdots\partial y_{j_\ell}}
       = \Dfrac{\partial^\ell f(T\yvec)}{\partial y_{j_1}\cdots\partial y_{j_\ell}}
      =
      \sum_{u_1,\ldots,u_\ell\in [n]}
         t_{u^{}_1j^{}_1}\!\cdots t_{u^{}_\ell j^{}_\ell} \,
         \Dfrac{\partial^\ell f(\xvec)}{\partial x_{u_1}\!\cdots\partial x_{u_\ell}}.
    \]
    Consequently,
  \[
     \widehat\Delta_\ell(g) \leq \max_{j\in[n]} 
        \sum_{\substack{\{ j^{}_2,\ldots,j^{}_\ell\}\subseteq [n]\\
         j\notin \{ j^{}_2,\ldots,j^{}_\ell\}}}\;
        \sum_{u^{}_1,\ldots,u^{}_\ell\in [n]} 
        \abs{t_{u^{}_1j} t_{u^{}_2j^{}_2}\cdots t_{u^{}_\ell j^{}_\ell}}
       \sup_{\xvec\in \varOmega}\;
        \biggl| \Dfrac{\partial^\ell f(\xvec)}{\partial x_{u^{}_1}\!\cdots\partial x_{u^{}_\ell}} \biggr|.
  \]
  Summing over $\{j_2,\ldots,j_\ell\}$ yields
  \begin{align*}
     \widehat\Delta_\ell(g) &\leq \frac{\infnorm{T}^{\ell-1}}{(\ell-1)!}\; 
     \max_{j\in[n]}  \sum_{u^{}_1,\ldots,u^{}_\ell\in [n]} \abs{t_{u^{}_1j}}
     \sup_{\xvec\in \varOmega}\;
        \biggl| \Dfrac{\partial^\ell f(\xvec)}{\partial x_{u^{}_1}\!\cdots\partial x_{u^{}_\ell}} \biggr| \\
       & \leq  \Dfrac{\infnorm{T}^{\ell-1}}{(\ell-1)!}\;
         \medbar\Delta_\ell(f, \varOmega)   
         \,\max_{j\in [n]} \sum_{u^{}_1\in[n]} \abs{t_{u^{}_1j}} 
      \leq  \Dfrac{\infnorm{T}^{\ell-1}\,\onenorm{T}}{(\ell-1)!}\;
         \,\medbar\Delta_\ell(f, \varOmega).
  \end{align*}
  This completes the proof.
  \end{proof}

 Finally, we are ready to prove the main result of this section.

\begin{thm}
\label{gauss2pt}
Let $\eps>0$ be constant.
 Let $A$ be an $n\times n$ positive-definite symmetric real matrix
 and let $T$ be a real matrix such that $T\trans\! AT=I_n$.
  Suppose the following assumptions hold for some $m=m(n) \in [n]$, measurable set $\Medbar\varOmega \subseteq \Reals^n$,
  measurable function $f: \Reals^n\to\Reals$, 
  and
  $\rho_1, \rho_2, \rho_3, \alpha, b, \hat b, C, \widehat C \ge 0$
  are functions of $n$.
    \begin{itemize}\itemsep=0pt
      \item[(i)] $U_n(\rho_1)
      \subseteq T^{-1}(\Medbar\varOmega)
      \subseteq U_n(\rho_2)
      \subseteq U_n(\rho_3)
    $, where  $\rho_3 \ge \rho_2 \geq \rho_1 \geq
    (\log n)^{1/2+\eps}$.
%
  \item[(ii)] For $\ell\in [m]$,
    \[
     \max_{j\in[n]} \sum_{V\in\binom{[n]}{\ell}:j\in V} \Delta_V(f(T\xvec),U_n(\rho_2))
     \le \alpha <\dfrac{1}{400},
    \]
    where $\Delta_V$
  is defined in \eqref{def:Delta}.
     \item[(iii)] $ 
     \max\{bm+b^2m^2,\hat bm+\hat b^2m^2\}\le n$.
     \item[(iv)]
\[
   \abs{f(\xvec)} \le \begin{cases}
       \,Ce^{b\xvec\trans\! A\xvec/n}, & \text{for $\xvec\in\Reals^n$}; \\[1ex]
       \,\widehat Ce^{\hat b\xvec\trans\! A\xvec/n}, & \text{for $\xvec\in T(U_n(\rho_3))$}.
   \end{cases}
\]
    \end{itemize}
    Then for sufficiently large $n$,
     we have
     \[
        \int_{\Medbar\varOmega} e^{-\xvec\trans\!A\xvec + f(\xvec)}\,d\xvec
        =
          \pi^{n/2}\abs{A}^{-1/2}
   \exp\bRg{ \,\sum_{r=1}^m \Dfrac{\kappa_r (f(\X))}{r!} + \delta }, 
     \]     
     where $\X$ is a random vector with the Gaussian density
    $\pi^{-n/2} \abs{A}^{1/2} e^{-\xvec\trans\!A\xvec}$, and
     \[
            \delta   \leq 
            2 n(100 \alpha)^{m+1} +e^{-\rho_1^2/2}
      +6^{m}
\bR{
(C^m + C )
e^{- (\rho_3^2-bm)/2 }
+ 
( \widehat C^m + \widehat C )
e^{-  (\rho_1^2-\hat bm)/2 }
}.
     \]
\end{thm}

\begin{proof}
  Let $\yvec = T^{-1}\xvec$.  Since $T\trans\! A T = I_n$, we have $\abs{T}=\abs{A}^{-1/2}$ and
  \begin{equation}\label{eq:rotate}
     \int_{\Medbar\varOmega} e^{-\xvec\trans\!A\xvec + f(\xvec)}\,d\xvec
     = \abs{A}^{-1/2} \int_{T^{-1}(\Medbar\varOmega)} e^{-\yvec\trans\yvec + f(T\yvec)}\,d\yvec.
  \end{equation}  
Let  $\rho \in \{\rho_1,\rho_2\}$.
Let $\Y$ have the Gaussian density $\pi^{-n/2}e^{-\yvec\trans\yvec}$ and define $p:=\P{\Y\notin U_n(\rho)}$.
Using cumulant expansion in Theorem \ref{cumuThm}, and noting that 
$(1-p)\pi^{n/2} 
= \int_{U_n(\rho)} e^{-\yvec\trans\yvec}\,d\yvec$,
we obtain
  \[
     \log \biggl( 
     \Dfrac{1}{(1-p)\pi^{n/2}} \int_{U_n(\rho)} e^{-\yvec\trans\yvec+f(T\yvec)}\,d\yvec
       \biggr)= 
        \sum_{r=1}^m 
        \Dfrac{1}{r!} 
      \kappa_r \bigl(f(T\Y) \mid \Y\in U_n(\rho) \bigr)
 + \delta',  
  \]
  where $\delta' = n \log (1+\delta)$
  for some $\abs\delta \leq e^{(100\alpha)^{m+1}}-1
  \le 1/2$.
Since $\abs{\log(1+x)}
\le 2\log(1+\abs{x})$
for $\abs{x} \le 1/2$,
we have
\[
   |\delta'|  =
   |n \log (1+\delta)| 
   \leq 2 n(100 \alpha)^{m+1}.
\]

Next, we relate the truncated cumulants
to cumulants, and bound the error.
By standard
  bounds on the tail of the Gaussian distribution,
  we have
  $
  p\le 2ne^{-\rho^2}
  $.
  Under our assumptions, for sufficiently large $n$,
  we have
  \[
   \abs{\log(1 - p)}
   \leq e^{-\rho^2/2}.
  \]

Now apply Lemma~\ref{truncated} to $f(\xvec)$, with $\varOmega=T(U_n(\rho))$,
and $\widehat\varOmega=T(U_n(\rho_3))$.  Note that 
$n\gg bm$, $p=\pr(\X\notin \varOmega)\le 2ne^{-\rho^2}$ and $\hat p=\pr(\X\notin\widehat\varOmega)\le 2ne^{-\rho_3^2}$.
This gives
\begin{align*}
\sum_{r=1}^m \Dfrac{1}{r!}\,
&
\left| \kappa_r \bigl(f(T\Y) \mid \Y\in U_n(\rho) \bigr) -\kappa_r (f(T\Y)) \right|  \\
&\le\sum_{r=1}^m
6^{r+1}
\bR{
C^r e^{br/2} \hat p^{1-br/n}
+ 
2\widehat C^r e^{\hat br/2} p^{1-\hat br/n} } \\
&\le m 6^{m+1}
\bR{
( C^m + C ) e^{bm/2} \hat p^{1-bm/n}
+ 
( \widehat C^m + \widehat C )
2 e^{\hat bm/2} p^{1-\hat bm/n} } \\
&\le
6^{m} \bR{ (C^m + C )
e^{- (\rho_3^2-bm)/2 }
+ 
( \widehat C^m + \widehat C )
e^{-  (\rho^2-\hat bm)/2 }
}.
\end{align*}

In view of the constraints on $\rho_1, \rho_2$,
we have
  \begin{equation}\label{newgen2}
     \int_{U_n(\rho)} e^{-\yvec\trans\yvec+f(T\yvec)}\,d\yvec = 
           \pi^{n/2} 
           \exp\bRg{\, \sum_{r=1}^m \Dfrac{\kappa_r (f(T\Y))}{r!} + \delta_\rho } ,
  \end{equation}
  where 
  \[
      \abs{\delta_\rho}       
      \le 2 n(100 \alpha)^{m+1} +e^{-\rho^2/2}
      +6^{m}
\bR{
(C^m + C )
e^{- (\rho_3^2-bm)/2 }
+ 
( \widehat C^m + \widehat C )
e^{-  (\rho^2-\hat bm)/2 }
}.
  \]
By assumption (i) and the fact that the integrand is positive, we obtain
\[
\int_{U_n(\rho_1)} e^{f(T\yvec)-\yvec\trans\yvec}\,d\yvec 
\leq 
 \int_{T^{-1}(\Medbar\varOmega)} e^{f(T\yvec)-\yvec\trans\yvec}\,d\yvec 
 \leq 
 \int_{U_n(\rho_2)} e^{f(T\yvec)-\yvec\trans\yvec}\,d\yvec .
\]
Using  \eqref{eq:rotate} and applying \eqref{newgen2} twice with $\rho =\rho_1$ and $\rho=\rho_2$,
the proof is complete.
\end{proof}


\section{Eulerian orientations
and the proof of Theorem~\ref{tm-eo}}\label{S:proof-main}

Throughout this section, we work under
the assumptions of Theorem \ref{tm-eo}.
For $\thetavec \in (\Rmodpi)^{n}$,
define 
\begin{equation}
F(\thetavec) := \prod_{jk \in G} \cos(\th_j - \th_k).
\label{Fdef}
\end{equation}
Recalling
\eqref{EOinitial2}, we have
\begin{equation}\label{NGint}
\EO(G) 
= 2^{ |E(G)| } \pi^{-n}   \int_{(\Rmodpi)^n} F(\thetavec)\,d\thetavec.
\end{equation}
The dominating contribution to \eqref{NGint} is 
from the region $\varOmega_0$ defined below,
where all $\theta_j$ are approximately the same, 
which makes it possible to estimate the integral asymptotically.

Given $x\in\Reals$, 
an \textit{interval} of $\Rmodpi$ of \textit{length} $\rho\ge 0$ is a set of the form
\begin{equation}
\label{Idef}
    I(x,\rho) := \{ \theta \in \Rmodpi \st \semiabs{x-\theta} \le \nfrac 12\rho \},
\text{~~~where~~~}
    \semiabs{x} := \min\{\abs{x-k\pi} \st k\in \Integers \}.
\end{equation}
Recall that $\dmax $ is the maximum degree of $G$.
Define
\begin{equation}\label{def:rho0}
\rho_{0} := 
\zeta^2 \dmax^{-1/2} \log^{3/2} n,
\end{equation} 
where
\begin{equation}\label{def:zeta}
     \zeta := \min\{ \dmax^{1/24}\log^{-1/3}n, (\log\log n)^{1/3}\}.
\end{equation}
Let  $\varOmega_0$ be the region consisting of those $\thetavec \in (\Rmodpi)^n$
such that all components $\theta_j$ can be covered by an interval of $\Rmodpi$ of length
at most  $\rho_{0}$:
\begin{equation}\label{def:omega0}
    \varOmega_0 := \bigl\{\thetavec \in (\Rmodpi)^n \st 
      \text{there exists $x\in\Rmodpi$ such that $\thetavec\in I(x, \rho_{0} )^n$}\bigr\}.
\end{equation}
Define
\[
     J_0 :=   \int_{  \varOmega_0} F(\thetavec)\,
             d\thetavec.
\]

The proof of Theorem~\ref{tm-eo} consists of two parts. In Section \ref{S:inside}, we estimate $J_0$ using Theorem \ref{gauss2pt} and some preliminary lemmas given in Section \ref{S:prelim}.
Then, in Section \ref{S:outside} we show that the integral over the region $(\Rmodpi)^n \setminus \varOmega_0$  is negligible with respect to $J_0$. Finally, we combine all the parts in Section \ref{S:proof-main-main} 
and also prove Corollary \ref{cor-eo}.

\subsection{Preliminaries}\label{S:prelim}

Here we state  some lemmas from \cite{isaev2020asymptotic,isaev2018complex} and also prove an additional lemma that is needed for the proof of Theorem \ref{gauss2pt}.
If $T:\Reals ^n\rightarrow \Reals ^n$ is a linear operator, let 
$\ker T = \{\x\in \Reals ^n \st T\x = 0\}$.

The first lemma helps to deal with integrals over a subspace of $\Reals^n$.
\begin{lemma}[{\cite[Lemma 4.6]{isaev2018complex}}]
\label{LemmaQW}
  Let $Q,W:\Reals ^n\rightarrow \Reals ^n$ be linear operators 
  such that $\ker Q \cap \ker W = \{\boldsymbol{0}\}$ and ${\rm span}(\ker Q, \ker W) = \Reals^n$. 
  Let $n_\perp$  denote the dimension of\/ $\ker Q$.
  Suppose $\varOmega \subseteq \Reals^n$ and
  $F:\varOmega\cap Q(\Reals^n) \to\Complexes$.
  For any $\eta>0$, define
  \[
   \varOmega_{\eta} = \bigl\lbrace \xvec\in\Reals^n \st
     Q\xvec\in \varOmega \text{~and~}
     W\xvec\in U_n(\eta) \bigr\rbrace.
  \]
   Then, if the integrals exist,
\begin{equation}
\int_{\varOmega \cap Q(\Reals^n)} F(\yvec)\,d\yvec
    = (1 - \delta)^{-1}\,\pi^{-\nperp/2} \,\Abs{Q\trans Q + W\trans W}^{1/2}
     \int_{\varOmega_{\eta}} F(Q\xvec)\, e^{-\xvec\trans W\trans W \xvec}
        \,d\xvec,
\label{eq:int-lm}
\end{equation}
   where 
\[
0\le \delta < \min(1,n e^{-\eta^2}). 
\]
      Moreover, if  $U_n(\eta_1) \subseteq \varOmega \subseteq  U_n(\eta_2)$ for some $\eta_2 \geq \eta_1 >0$ then
  \[
      U_n\biggl(\min\Bigl(\Dfrac{\eta_1}{\infnorm{Q}},
          \Dfrac{\eta}{\infnorm{W}} \Bigr)\biggr)
        \subseteq \varOmega_{\eta} \subseteq
      U_n\Bigl( \infnorm{P}\, \eta_2 +  \infnorm{R}\, \eta \Bigr)
  \]	
   for any linear operators $P,R : \Reals^n \rightarrow \Reals^n$ such that $PQ + RW$ is equal to the identity operator on $\Reals^n$.     
\end{lemma}

Let 
\begin{equation}
A := \dfrac{\dmax}{n} J_n + \dfrac12 L,
\label{afdefs}
\end{equation}
 where $L=L(G)$ is the Laplacian matrix of $G$ defined by \eqref{Ldef},
and $J_n$ denotes the $n\times n$ matrix with every entry one.
The next lemma will be useful for the norm bounds required in the application of Theorem \ref{gauss2pt}.
\begin{lemma}[{\cite[Lemma 12]{isaev2020asymptotic}}]
\label{l:norms}
Under the assumptions of Theorem \ref{tm-eo}, the following estimates hold.  
\begin{itemize}\itemsep=0pt
  \item[(a)] $\infnorm{A^{-1}}=O\(\dmax^{-1}\lognd \)$.
  \item[(b)] If $A^{-1}=(a_{jk})$, then $a_{jj}=O(\dmax^{-1})$ and
    $a_{jk}=O\(\dmax^{-2}\lognd \)$ uniformly for  $j\ne k$.
  \item[(c)] There exists a symmetric positive-definite matrix $T = A^{-1/2}$ such that
    $T\trans\!AT=I_n$.  Moreover, $\infnorm{T}=O(\dmax^{-1/2}\log^{1/2}\!n)$
    and $\infnorm{T^{-1}}=O(\dmax^{1/2})$.
\end{itemize}
\end{lemma}

The next three lemmas are needed for estimates outside the region $\varOmega_0$.

\begin{lemma}[\cite{isaev2020asymptotic}]
\label{separator}
 Suppose $0<t<\frac 13\pi$ and $q\le \frac 15n$.
 Let $X=\{x_1,\ldots,x_n\}$ be a multisubset of $\Rmodpi$ such that
 no interval of
 length $3t$ contains $n-q$ or more elements of~$X$.
 Then there is some interval $I(x,\rho)$, $\rho<\frac13\pi$, such that both 
 $I(x,\rho)$ and $\Rmodpi-I(x,\rho+t)$ contain at least $q$ elements
 of~$X$.
\end{lemma}

\begin{lemma}[\cite{isaev2020asymptotic}]
\label{l:paths} Let $G$ be a graph with maximum degree $ \dmax $. Assume also that
$h(G)\ge \gamma  \dmax $ for some $\gamma>0$.  
Let $U_1, U_2$ be two disjoint sets of vertices.
Then, there exist at least 
\[
\gamma  \dmax \, \Dfrac{ \min\{|U_1|,|U_2|\}}{2\ell(U_1,U_2)}
\]
 pairwise edge-disjoint paths in $G$
  with  one end in~$U_1$ and the other end in $U_2$,
  and of lengths  
  bounded above by 
\[
	\ell(U_1,U_2) := 2+2\log_{1+\gamma/2} 
	  \biggl(\Dfrac{|V(G)|}{  \min\{|U_1|,|U_2|\}+  \gamma   \dmax /2}\biggr).
\]  
  
\end{lemma}

\begin{lemma}\label{l:fapp}
 $\abs{\cos x }$ is a decreasing function of $\semiabs{x}$
   and
\begin{equation}\label{fapprox}
  \cos^2 x
   \leq 
                   \exp\bigl( - \dfrac{1}{4} \semiabs{x}^2 \bigr). 
\end{equation}
In addition, for any $\semiabs{y} \le \semiabs{x}$, we have 
\begin{equation}\label{fapprox2}
	\abs{\cos x}
	\le \abs{\cos y}\, \Exp{- \dfrac{1}{4\pi} (\semiabs{x}^2 - \semiabs{y}^2) (\pi - \semiabs{x} - \semiabs{y}))}.
\end{equation}
\end{lemma}
\begin{proof}
The first part of~\eqref{fapprox} follows from the definition of $\cos(x)$
and observing that
$\abs{\cos(x)}=\cos(\semiabs{x})$ for all $x$.
Therefore we can assume that $0\le y\le x\le \frac12\pi$,
which implies that $\semiabs{x}=x$ and $\semiabs{y}=y$.

By the concavity of $\cos z$ on $[0,\frac{\pi}{2}]$, we have
$\cos z \geq 1 - \frac{2z}{\pi}$ in this range, which in turn implies
(by symmetry about the line $z = \frac{\pi}{2}$) that  
\begin{equation} \label{sinlower}
    \sin{z} \geq \dfrac1\pi z(\pi-z), \quad z \in [0, \pi ].
\end{equation}
Therefore, for $0\le x\le \frac12\pi$,
we have that
\[
\cos^2 x
= 1 - \sin^2 x
\le 1 - \dfrac{1}{\pi^{2}} x^{2}(\pi-x)^{2}
\le 1 - \dfrac{1}{4} x^{2}
\le \exp\Bigl(-\Dfrac{x^2}{4}\Bigr).
\]

Inequality~\eqref{fapprox2} is trivial if $x=y$, so assume that
$0\le y < x\le\frac12\pi$. 
In that case, $\cos(y)\ne 0$ and
\[
      \Dfrac{ \cos^{2}x}{ \cos^{2}y }
      \le  \Dfrac{1- \sin^2 x}{1-\sin^2 y}
      \le \exp \biggl( -\Dfrac{\sin^2 x-\sin^2 y}{1-\sin^2 y} \biggr)
      \le \exp \( -(\sin^2 x-\sin^2 y) \).
\]
Finally, by \eqref{sinlower} with $z= x+y$ and $z=x-y$, we have
\begin{align*}
  \sin^2 x - \sin^2 y 
  = \sin{(x+y)}\sin{(x-y)} 
  &\geq \dfrac{1}{\pi^2} (x^2 - y^2) (\pi - x + y) (\pi - x - y) \\
  &\geq \dfrac{1}{2\pi} (x^2 - y^2) (\pi - x - y)
\end{align*}
for $0\le y\le x\le \frac12\pi$, which completes the proof of~\eqref{fapprox2}.
\end{proof}

\subsection{The integral inside $\varOmega_0$}\label{S:inside}

By the definition of $\varOmega_0$ in \eqref{def:omega0}, 
we have that
\[
    \{ \thetavec\in (\Rmodpi)^n\st \max_{j\in [n]}\,\semiabs{\theta_j-{}\theta_n}
        \le \rho_{0}   \}  
    \subseteq \varOmega_0\subseteq
     \{ \thetavec\in (\Rmodpi)^n \st 
     \max_{j \in [n]}\, \semiabs{\theta_j-\theta_n}
        \le 2 \rho_{0}\}.
\]
Observe from the definition of $F(\thetavec)$
in \eqref{Fdef} that, for any $\thetavec \in \Reals^n$,
\[
   F(\thetavec) = F(\thetavec - \theta_n \boldsymbol{1}),
\]
where $\boldsymbol 1:=(1,\ldots,1)\trans  \in \Reals^n$.
Therefore, 
\[
   J_0=\pi \int_{\varOmega_0 \cap \mathcal{L}} F(\thetavec')\, d\thetapvec,
\]
where $\mathcal{L}:= \{\thetavec' \in \Reals^n \st \theta_n' =0\}$ and
\[
U_{n-1}(\rho_{0})
\subseteq
( \varOmega_0 \cap \mathcal{L} )
\subseteq U_{n-1}(2\rho_{0}).
\]

Next, we lift the integral back to full dimension using Lemma \ref{LemmaQW}.
Let $S$ be the matrix with ones in the last column and zero elsewhere.
Define
\[
   P = I_n - \dfrac1n J_n,~~
   Q = I_n - S,~~
   R = \dmax^{-1/2} I_n, ~~\text{ and }~~
   W = \dmax^{1/2}n^{-1} J_n.
\]
One can easily check that $PQ+RW=I$, and also that
$\ker Q\cap\ker W=\{\boldsymbol0\}$, $\ker Q$ has dimension~1 and
$\sp(\ker Q,\ker W)=\Reals^n$.
From  \cite[Section 3.1]{isaev2020asymptotic}, we know that
\[
\abs{Q\trans Q+W\trans W}=n\dmax,  
 \quad \infnorm{P}\le 2\quad \infnorm{Q} = 2, \quad 
\infnorm{R}=\dmax^{-1/2}, \quad \infnorm{W}=\dmax^{1/2}.
\]
Then applying Lemma \ref{LemmaQW} with
$\eta_1 =\rho_0$ and $\eta_2=2\rho_0$, we find, from \eqref{eq:int-lm}, that
\begin{equation}\label{eq:J0_2}
   J_0 = \(1+e^{-\omega(\log n)}\)\, \pi^{1/2} (\dmax n)^{1/2}
            \int_{\varOmega}
             \widehat F(\thetavec)\,
             d\thetavec, 
   \end{equation}          
    where $\varOmega$ is  some  region such that  
    \begin{equation}
    U_n\left( \rho_{0} /2  \right) \subseteq \varOmega
    \subseteq U_n(5 \rho_{0})
    \label{EO-varOmega}
    \end{equation}
          and
    \[ 
        \widehat F(\thetavec) 
        := \exp\biggl( -\Dfrac \dmax n
        \Bigl( \sum_{i \in [n]} \theta_i \Bigr)^2
        + \sum_{jk \in G} \log \cos(\th_j - \th_k)
        \biggr).
\]

Now we are ready to derive the main estimate of this section.
 \begin{thm}\label{l:inside}
 Under the assumptions of Theorem \ref{tm-eo}, we have
\[
    J_0 = \pi^{(n+1)/2} \dmax^{1/2}
               n^{1/2} \abs{A}^{-1/2}
              \exp\bRgg{ \,
\sum_{r=1}^{M}  \Dfrac{1}{r!}\, \kappa_r (f_{K}(\X))
+ \O{ n^{-c} } }.
\]
\end{thm}

\begin{proof}
Define
\[
K_0 := \biggl\lceil\frac{2(c+1)\log n}{\log \dmax}\biggr\rceil.
\]
In view of the definition
of $f_{K_0}$ in \eqref{eq:eoTaylor}, 
and the definition of matrix $A$ in \eqref{afdefs}, 
using Taylor's series, we find that, for all  $\thetavec \in \varOmega$,
\begin{equation}\label{log-expansion}
\begin{aligned}
\log 
\widehat F(\thetavec) 
&= -\dfrac \dmax n
        \Bigl( \sum_{i \in [n]} \theta_i \Bigr)^2
-  \dfrac{1}{2} \sum_{jk \in G} (\theta_j-\theta_k)^{2}
+ \sum_{t\ge2} c_{2t} \sum_{jk \in G} (\theta_j-\theta_k)^{2t} \\ 
&=-\thetavec \trans\!A\thetavec  + f_{K_0}(\thetavec) + O(n^{-c}).
\end{aligned}
\end{equation} 
In the above estimation,
we recall from \eqref{def:rho0}  that $\rho_0 = 
\zeta^2 \dmax^{-1/2} \log^{3/2}n$
and observe that the bound on $c_{2t}$ in \eqref{c-bound} 
implies
\[
\sum_{t > K_0} c_{2t} \sum_{jk \in G} (\theta_j-\theta_k)^{2t}
=   
\sum_{jk \in G} O( (2\rho_0)^{2(K_0+1)}) 
=
O\(n\dmax \( 2\zeta^2 \dmax^{-1/2} \log^{3/2} n\)^{2K_0+2}\)
=O(n^{-c}).
\]

Then the argument continues by combining \eqref{log-expansion} and   Theorem \ref{gauss2pt}. 
First, we verify its assumptions.
By Lemma~\ref{l:norms}(c), 
there is a symmetric positive-definite matrix $T$ such that $T\trans \!AT = I_n$ and 
\begin{equation}\label{proof:normT}
        \norm{T}_\infty = O(\dmax^{-1/2}\log^{1/2}n),\qquad 
        \norm{T^{-1}}_\infty = O(\dmax^{1/2}).
\end{equation}
 Then, combining  \eqref{EO-varOmega} and \eqref{proof:normT}, we obtain 
\[
 U_n(\rho_1) \subseteq T^{-1}(\varOmega)\subseteq U_n(\rho_2), 
\]
for some $\rho_1,\rho_2$ such that
\[
\rho_1 = \Theta\( \zeta^2 \log n\) 
\quad\text{ and }\quad
\rho_2 = \Theta\( \zeta^2 \log^{3/2}n\).
\]
This is because for $\thetavec' \in U_n( \O{ \rho_{0} } )$, we have that
\[
\infnorm{T^{-1}\thetavec'}
\le \infnorm{T^{-1}}\, \infnorm{\thetavec'}
= O\( \zeta^2 \log^{3/2} n\),
\]
and
\[
\infnorm{T}^{-1}\,\infnorm{\thetavec'}
= 
\O{ \Dfrac{ \zeta^2 \dmax^{-1/2} \log^{3/2} n }
{\dmax^{-1/2}\log^{1/2}\!n} } = O( \zeta^2 \log n).
\]

Next, we estimate
$\medbar{\Delta}_\ell(f_{K_0}, T(U_n(\rho_2)))$ for $\ell \ge 1$.
Using \eqref{proof:normT}, we find that, for any $\thetavec \in T(U_n(\rho_2))$,
\[
    \|\thetavec\|_\infty = O\(\rho_2 \dmax^{-1/2} \log^{1/2}n \) 
    = O\( \zeta^2 \dmax^{-1/2} \log^{2} n\)
    \ll K_0^{-1}.
\]
Then, noting from \eqref{c-bound}
that
$c_{2t}<1$,
for any integer $i\ge 1$,
\begin{align*}
\Dfrac{\partial^{i} f_{K_0}(\thetavec)}{\partial^{i} \theta_j}
&= \sum_{t=2}^{K_0} ( 2t )_i\, c_{2t} \sum_{k: jk \in G} (\theta_j-\theta_k)^{ (2t-i)_{+} }
= O\( 
i!\,\dmax\, \infnorm{\thetavec}^{ (4-i)_{+} } \),
\end{align*}
where $(\cdot)_i$ denotes the falling factorial and $(x)_{+}$ denotes $\max \{x, 0\}$.
Similarly, for $jk \in G$ and distinct $p, q \in [K_0]$,
\begin{align*}
\Dfrac{\partial^{p+q} f_{K_0}(\thetavec)}{\partial^p \theta_j\partial^q \theta_k}
&= \sum_{t=2}^{K_0}\, (2t)_{p+q} c_{2t} (\theta_j-\theta_k)^{(2t-p-q)_{+}}
= O\( (p+q)!\,
\infnorm{\thetavec}^{ (4-p-q)_{+} } \).
\end{align*}
All higher-order mixed derivatives with three or more distinct indices are zeros.
Since $T$ is symmetric, we have $ \onenorm{T} = \infnorm{T}$. 
Then, using \eqref{proof:normT}, we obtain, uniformly for $\ell\ge 1$,
\begin{align*}
& (2\rho_2)^\ell \Dfrac{\infnorm{T}^{\ell-1}\, \onenorm{T}}
 {(\ell-1)!}
\max_{u^{}_1\in[n]} 
\sum_{u^{}_2,\ldots,u^{}_\ell\in [n]}\,
          \sup_{\thetavec\in T(U_n(\rho_2))}\;
        \biggl| \Dfrac {\partial^\ell f_{K_0}(\thetavec)}{\prod_{r\in [\ell]}\partial \theta_{u^{}_r}}\biggr| \\
&=  
\Dfrac{\(
2\rho_2 \norm{T}_\infty \)^{\ell}}
 {(\ell-1)!}
\max_{u^{}_1\in[n]} 
\biggl(\,
         \sum_{u^{}_2 = u_1 }
         \, \sup_{\thetavec\in T(U_n(\rho_2))}\;
        \biggl| \Dfrac {\partial^\ell f_{K_0}(\thetavec)}{\prod_{r\in [\ell]}\partial \theta_{u^{}_r}}\biggr| +\sum_{u^{}_2 \in [n]\setminus\{u_1\} }
          \sup_{\thetavec\in T(U_n(\rho_2))}\;
        \biggl| \Dfrac {\partial^\ell f_{K_0}(\thetavec)}{\prod_{r\in [\ell]}\partial \theta_{u^{}_r}}\biggr| 
        \,\biggr) \\
&=  
O\( \(
2\rho_2 \norm{T}_\infty \)^{\ell}
\ell
\dmax \infnorm{\thetavec}^{ (4-\ell)_{+} } \)
= O\(\dmax^{-1} \zeta^8 \log^{8} n\)
= o(1),
\end{align*}
in view of the definition of $\zeta$ in \eqref{def:zeta}.
Then by Lemma \ref{Deltabound}, we have
\[
     \max_{j\in[n]} \sum_{V\in\binom{[n]}{\ell}:j\in V} \Delta_V(f(T\xvec),U_n(\rho_2))
     \le \alpha <\dfrac{1}{400},
\]
where $100\alpha = e^{\beta} \zeta^8 \dmax^{-1} \log^{8}n
$ for some constant~$\beta$.

Next,
using Lemma \ref{l:norms}(a), we find that 
\[
\norm{\thetavec}_2^2
\le
\thetavec\trans  \!A \thetavec 
\norm{A^{-1}}_\infty
=
O\Bigl(
\thetavec\trans  \!A \thetavec \,
\dmax^{-1}\lognd
\Bigr).
\]
Then for all $\thetavec\in\Reals^n$,
we have
\begin{align*}
\abs{f_{K_0}(\thetavec)} 
&= 
\biggl|\,\sum_{t=2}^{K_0} c_{2t} \sum_{jk \in G} (\theta_j-\theta_k)^{2t}\,\biggr| \\
&\le 
\O{n\dmax K_0} 4^{K_0}
\bR{
\norm{\thetavec}_2^4
+ \norm{\thetavec}_2^{2K_0} } \\
&
\le 
\O{n\dmax K_0} 4^{K_0}
\bR{ 
(\thetavec\trans  \!A \thetavec \norm{A^{-1}}_\infty )^2
+
(\thetavec\trans  \!A \thetavec \norm{A^{-1}}_\infty )^{K_0} } \\
&
\le \O{n\dmax K_0} 4^{K_0}
\bR{ 
\bR{ O(1)\thetavec\trans  \!A \thetavec 
\dmax^{-1} \lognd }^2
+ 
\bR{ O(1)\thetavec\trans  \!A \thetavec \dmax^{-1} \lognd }^{K_0} } \\
&
\le 
o(n^{K_0})
\bR{ 
\dfrac{1}{2}
\bR{ \dfrac1n
\thetavec\trans \!A \thetavec  }^2
+ 
\dfrac{1}{K_0!}
\bR{ 
\dfrac1n\thetavec\trans  \!A \thetavec
}^{K_0} } \\
&\le 
n^{\O{\log n/\log \dmax}}
 e^{\frac{1}{n}
 \thetavec\trans  \!A \thetavec}.
\end{align*}
If $\infnorm{T^{-1}\thetavec} \le \dmax^{1/2}\log^{-1/2}n$, then
$\infnorm{\thetavec} = O(1)$.
In that case, $\abs{f_{K_0}(\thetavec)}=O(n^3)$.
Let
\[
    M_0 :=
\biggl\lceil \Dfrac{ 2(c+1)\log n }{\log \dmax
- 8 \log\log n } \biggr\rceil.
\]
Now apply Theorem \ref{gauss2pt} with
\begin{align*}
\rho_3&=\dmax^{1/2}\log^{-1/2}n, 
&
\Medbar\varOmega&=\varOmega,
& m&=M_0,\\
C&=n^{O(\log n/\log \dmax)},
\widehat C=O(n^3),
&
b&=1,
\hat b=0,
&
\alpha&= \dfrac{1}{100}e^{\beta}\zeta^8 \dmax^{-1} \log^8 n.
\end{align*}
Combined with \eqref{eq:J0_2} and \eqref{log-expansion},
and using $\zeta\le \dmax^{1/24}\log^{-1/3}n$
from \eqref{def:zeta} give
\[
 J_0 = \pi^{(n+1)/2} \dmax^{1/2}
               n^{1/2} \abs{A}^{-1/2}
              \exp\bRgg{\, 
\sum_{ r = 1 }^{M_0}  \Dfrac{1}{r!}\, \kappa_r \( f_{K_0}(\X) \) +\delta},
\]
where
\begin{align*}
\abs{\delta} 
&\le \Exp{ 2 n(100 \alpha)^{M_0+1} } - 1 \\ 
&{\quad}+ e^{-\rho_1^2/2}
      +6^{M_0}
\bigl(
(C^{M_0} + C )
e^{- (\rho_3^2-bM_0)/2 }
+ 
( \widehat C^{M_0} + \widehat C )
e^{-  (\rho_1^2-\hat bM_0)/2 }
\bigr) \\
&= O( n^{-c} ).
\end{align*}
Therefore,
\[
    J_0 = \pi^{(n+1)/2} \dmax^{1/2} n^{1/2} \abs{A}^{-1/2}
              \exp\bRgg{ \,
\sum_{r=1}^{M_0}  \Dfrac{1}{r!}\, \kappa_r (f_{K_0}(\X))
+ \O{ n^{-c} } }.
\]

Finally, we show that $M_0$ and $K_0$ can be reduced to $M$ and $K$, respectively.
From Theorem~\ref{tm-eo} assumption (A2), the minimum degree
of~$G$ is at least $\gamma \dmax$. Applying Theorem~\ref{T:cumulant_bounds}(b)
and Lemma~\ref{l:norms}(a), there is some constant $a_1$ such
that
\begin{equation}\label{eq:rcumbnd}
\sum_{r=M+1}^{M_0}
\Dfrac{1}{r!}\, \kappa_r(f_{K_0}(\X))
    \leq 
\sum_{r=M+1}^{M_0} 
a_1^r \Dfrac{(4r-1)!!}{r!}\, n \dmax^{-r} \log^{r-1}\!\dfrac{2n}{\dmax}.
\end{equation}
In this range, $\dmax\gg r\log n$, so the term $r=M+1$ dominates
within a constant.
Since $(4M+3)!!/(M+1)! = O(6^M M^M)$
and $\log M = \log\log n - \log\log \dmax + O(1)$,
we find that this term is $O(n^{-c})$.
Thus we may replace $M_0$ by $M$ with a loss of~$O(n^{-c})$.

Using Theorem~\ref{T:cumulant_bounds}(c) in similar manner, 
there is a constant $a_2$ such that
\begin{equation}\label{eq:rcumbnd2}
   \sum_{r=1}^{M}
   \Dfrac{1}{r!}\, \bigl| \kappa_r(f_{K_0}(\X)) - \kappa_r(f_{K+1}(\X))\bigr|
    \le
    \sum_{r=1}^{M}\Dfrac{(2K+4r-3)!!}{(r-1)!}a_2^r\, 
        n \dmax^{-r-K+1} \log^{r-1}\!\dfrac{2n}{\dmax}.
\end{equation}
Since $\dmax\gg (K+r)\log n$ in this range, the value of~\eqref{eq:rcumbnd2}
is dominated within a constant by the term $r=1$.
Since $(2K+1)!! 
= O( K^{K})$,
and $\log K = \log\log n - \log\log \dmax + O(1)$,
we find that this term is $O(n^{-c})$.
Thus we may replace $K_0$ by $K$ with a loss of~$O(n^{-c})$.
This completes the proof.
\end{proof}

In the next section, we will need the following crude bound on $J_0$, which is a simple consequence of Theorem \ref{l:inside}.

\begin{lemma}\label{J0rough}
Under the assumptions of Theorem \ref{tm-eo}, we have 
$
     J_0 = e^{O(n\log n)}. 
$
\end{lemma}

\begin{proof}[Proof of Lemma \ref{J0rough}]
Note that all of the eigenvalues of $A^{-1}$ are bounded below by $\infnorm{A}^{-1}$
and bounded above by $\infnorm{A^{-1}}$.
Using Lemma \ref{l:norms}, we have
$\abs{
A
}^{-1/2} = e^{\O{n \log n}}$.
Using \eqref{eq:rcumbnd}, \eqref{eq:rcumbnd2} and recalling $M= O(\log n)$, we find that all cumulants in the formula of Theorem \ref{l:inside} are $o(n)$. The claim follows.
\end{proof}

\subsection{The integral outside $\varOmega_0$}\label{S:outside}

It remains to show that the integral of $F(\thetavec)$ is negligible outside $\varOmega_0$.
We follow the approach of \cite[Section~3.3]{isaev2020asymptotic}
with slight modifications and improvements. 
Define
\[
	 \srho := \zeta
  \dmax ^{-1/2} \log^{1/2} n.
\]
Then
\[
\Dfrac{\brho}{\srho}
= \zeta \log n.
\]

First, we bound the integral of $\abs{F(\thetavec)}$ in   region
\[
	\varOmega_1 :=  \bigl\{ \thetavec \in (\Rmodpi)^n \st \text{for every
	   $\xi\in\Rmodpi$, we have  $\card{\{ j \st \theta_j\in I(\xi,\srho) \}} <\tfrac45n$} \bigr\},
\]
where $I(\cdot, \cdot)$ is defined by \eqref{Idef}.

\begin{lemma}\label{S1bound}
We have
\[
    \int_{\varOmega_1} \abs{F(\thetavec)}\,d\thetavec 
    = e^{- \omega(n\log n)} J_0.
\]
\end{lemma}

We need the following lemmas to prove 
Lemma \ref{S1bound}.

\begin{proof}[Proof of Lemma \ref{S1bound}]
   If $\thetavec\in \varOmega_1$, the definition of $\varOmega_1 $ implies
   that every interval of $\Rmodpi$ of
  length $\srho$ has fewer than $\frac45n$ components of~$\thetavec$.
  Applying Lemma~\ref{separator} with $t=\frac13 \srho$,
  $q=\frac15n$, and $X = \thetavec$ tells us that there exist $p\in \Rmodpi$ and $s<\frac\pi 3$
  such that both  $I(p,s)$ and $\Rmodpi-I(p,s+t)$ contain at least $\frac15n$
  components  of $\thetavec$.

  For such~$\thetavec$, 
let  $U, U'$ denote
 the indices of the elements of $\thetavec$ belonging to $I(p,s)$ and $\Rmodpi-I(p,s+t)$ respectively.
Then $\semiabs{\theta_j-\theta_k}\ge t$
  whenever $j\in U,k\in U'$. Note that $t=o(1)$.
Consider any of  the paths $v_0,v_1,\ldots,v_\ell$ provided by 
  Lemma~\ref{l:paths}.  
  Note that we have $\ell = \ell(U_1,U_2) = \O{1}$.
  By assumption,
  $\semiabs{\theta_{v_0}-\theta_{v_\ell}}\ge t$.
   Since
  $\semiabs{\,\cdot\,}$ is a seminorm, we have
  \[
  	\sum_{j=1}^\ell \,\semiabs{\theta_{v_j}- \theta_{v_{j-1}}}^2 
  	\ge   \Dfrac1\ell \biggl(\sum_{j=1}^\ell \semiabs{\theta_{v_j}- \theta_{v_{j-1}}}\biggr)^{\!2}  
	\ge  \Dfrac{t^{2}}{\ell}.
  \]
  Multiplying the bound~\eqref{fapprox} over all the edges of all the paths
  given by Lemma~\ref{l:paths} gives that
  \[
  |F(\thetavec)|
\le \Exp{
- \Omega\left( \gamma  \dmax \, \Dfrac{ \min\{|U_1|,|U_2|\}}{2\ell}\,
 \Dfrac{t^{2}}{\ell} \right)}
= \Exp{ - \Omega \left(  \Dfrac{ t^{2}  \dmax  n }{\ell^{2}}  \right) }.
  \]
Recalling 
$t=\frac13 \srho$
and 
$\srho
= \zeta \dmax ^{-1/2} \log^{1/2} n$,
using $\pi^n$ to bound the
  volume of $\varOmega_1$ yields
\begin{align*}
\int_{\varOmega_1} \abs{F(\thetavec)}\,d\thetavec 
&\le \pi^n \Exp{ - \Omega (  t^{2}  \dmax  n ) } \\
&= 
\pi^n \Exp{ - \Omega \left(  
\zeta^2 n \log n   \right) }
\le \pi^n \Exp{ - \omega \left(
n \log n \right) }
= e^{ - \omega \left(  n \log n \right) }.
\end{align*}  
  Then the result follows from Lemma~\ref{J0rough}.
\end{proof}

For  disjoint $U, W\subseteq V(G)$ 
define by $\varOmega_{U,W}$
the set of $\thetavec \in  (\Rmodpi)^n$, for which there exist some $x\in\Rmodpi$ and
$\rho$  with $\srho\le \rho \le \brho$ such that the following hold:
\begin{itemize}\itemsep=0pt
	\item[(i)] $\theta_j\in I(x,\srho)$ for at least $4n/5$ components $\theta_j$.
 	\item[(ii)]  $\theta_j \in I(x,\rho+\srho)$ if and only if $j \notin U$.
 	\item[(iii)] $\theta_j\in I(x,\rho+\srho)-I(x,\rho)$ if and only if $j \in W$.
 \end{itemize}

\begin{lemma}\label{l:UWcover}
We have
\[
	 (\Rmodpi)^n - \varOmega_0 - \varOmega_1 \subseteq \bigcup_{
	\substack{ 1\le |U|\le n/5 \\
	|W|\le 
    \zeta^{-1}|U| } } \varOmega_{U,W}.
\]
\end{lemma}
\begin{proof}
Any $\thetavec \in (\Rmodpi)^n -  \varOmega_1$ is such that 
  at least $4n/5$ of  its components  $\theta_j$  lie in some interval $I(x,\srho)$.
   Suppose it is not covered by any $\varOmega_{U,W}$.
   For $1 \leq k \leq \brho/\srho
   = \zeta \log n$, take $\rho = k\srho \le  \brho$ and let $U$
   correspond to the components not in $I(x, \rho +\srho)$.
Since  (iii) cannot hold, we obtain
\[
	\Dfrac{|\{j\st \theta_j \notin I(x,k\srho)\}|} {|\{j\st \theta_j \notin I(x, (k+1)\srho)\}|}
	= 1 + \Dfrac{|\{j \st \theta_j\in I(x,\rho+\srho)-I(x,\rho)\} |}
	{|\{j\st\theta_j \notin I(x,\rho+\srho)\}|} > 1 +
\zeta^{-1}.
\]
 Recalling that $|\{j\st\theta_j \notin I(x,\srho)\}| \le n/5$, we can apply this ratio repeatedly starting with $k=1$ to find that 
\[
	\Card{\{j\st\theta_j \notin I(x, \brho)\}} 
	  \le  \dfrac15 n \(1 + \zeta^{-1} \)^{- 
	   \zeta\log n+1}<1.
\]
This implies that $\thetavec \in \varOmega_0$, which completes the proof.
\end{proof}

Next, we bound the contribution of $\varOmega_{U,W}$ to the integral of \eqref{NGint} except the case when most of the coordinates  $\theta_j$ are highly concentrated.
 
\begin{lemma} 
\label{JUWbound}
	For any disjoint $U, W\subseteq V(G)$ with $|U|\le n/5$ and $|W|\le 
 \zeta^{-1} |U|$, we have
	\[
		\int_{\varOmega_{U,W}-\varOmega_2} 
		   \abs{F(\thetavec)}\, d\thetavec = e^{-\omega(|U| \log n)} 
		J_0,
	\]
 where
 \[
   \varOmega_2 :=  \bigl\{ \thetavec \in (\Rmodpi)^n \st \text{there is~}
         x\in\Rmodpi \text{ such that }
       \Card{\{ j \st \theta_j\in I(x,e^{- 
    \zeta \log n}) \}} \ge\tfrac45n \bigr\}. 
\]
\end{lemma}

\begin{proof}
Let $X:=V(G)-(U\cup W)$ and define the map 
$\phivec =(\phi_1,\ldots,\phi_n): \varOmega_{U,W} \to \varOmega_0$  as follows.
By the definition of $\varOmega_{U,W}$, for any $\thetavec\in\varOmega_{U,W}$
there is some interval of length at most $\brho$ that contains $\{\theta_j\}_{j\in X}$.
Let $I(z,\xi)$ be the unique shortest such interval. We can ignore parts of
$\varOmega_{U,W}$ that lie in $\varOmega_2$,
which means that we can assume $\xi\ge 
e^{- \zeta \log n}$.

Identifying $\Rmodpi$ with $(z-\frac12\xi,z-\frac12\xi+\pi]$, define
\[
  \phi_j  = \phi_j(\thetavec) :=
    \begin{cases}
	z+\dfrac12\xi - \dfrac{\xi}{\pi-\xi}\(\theta_j-z-\dfrac12\xi\),
			       & \text{ if }  j \in U \cup W;\\
	\theta_j, &  \text{ if } j \in X.
    \end{cases}
\]
For $j\in U\cup W$, $\theta_j\notin I(z,\xi)$ and $\phi_j$ maps the
complementary interval $I(z+\frac12\pi,\pi-\xi)$ linearly onto $I(z,\xi)$
(reversing and contracting with $z-\frac12\xi$
and
$z+\frac12\xi$ fixed).
For $j\in X$, $\theta_j\in I(z,\xi)$ and $\phi_j=\theta_j$.

Thus $\semiabs{\phi_j - \phi_k} \le\semiabs{\theta_j - \theta_k}$ 
for all $j,k$. From Lemma~\ref{l:fapp}, we find that 
\[	
	|\cos(\theta_j - \theta_k)|
\le |\cos(\phi_j - \phi_k)|.
\]
Moreover, for $j\in U$ and $k\in X$, we have
$\semiabs{\phi_j - \phi_k} \le  \semiabs{\theta_j - \theta_k} -\frac12\srho$.
Observing also that $\semiabs{\phi_j - \phi_k} \le \xi = o(1)$ 
and  using \eqref{fapprox2}, we find that
\[
	\Dfrac
	{|\cos(\theta_j - \theta_k)|}
	{|\cos(\phi_j - \phi_k)|}
	 \le e^{-\Omega(\srho^2)}
	 = e^{-\Omega( \zeta \dmax ^{-1} \log n )}.
\]
By assumption (2) of Theorem~\ref{tm-eo}
and the definition of the Cheeger constant, 
this bound applies to at least 
\[
h(G)|U| -  \dmax  |W| 
\ge (\gamma+o(1))  \dmax  |U|
\] 
pairs  $jk\in\partial_G U$, thus 
\[
	|F(\thetavec)| 
	= \biggl| \prod_{jk \in G} \cos(\th_j - \th_k) \biggr|
	= e^{-\Omega( \zeta^{2} |U| \log n)} |F(\phivec(\thetavec))|.
\]

Note that the map $\phivec$ is injective, since $I(z,\xi)$ can be determined
from $\{\phi_j\}_{j\in X}=\{\theta_j\}_{j\in X}$.
Also, $\phivec$  is analytic except at places where the map from
$\{\theta_j\}_{j\in X}$ to $(z,\xi)$ is non-analytic, which happens
only when two distinct components $\theta_j,\theta_{j'}$ for $j,j'\in X$
lie at the same endpoint of $I(z,\xi)$.
Thus, the points of non-analyticity of $\phivec$ lie on a finite number of hyperplanes,
which contribute nothing to the integral.
To complete the calculation, we need to bound the Jacobian of the transformation
$\phivec$ in the interior of a domain of analyticity.

We have
\[
    \frac{\partial \phi_j}{\partial\theta_k} =
      \begin{cases}
         1, & \text{~ if $j=k\in X$}; \\
          \pm\dfrac{\xi}{\pi-\xi}, & \text{~ if $j=k\notin X$}; \\[-0.5ex]
         0, & \text{~ if $j\ne k$ and either $j\in X$ or $k\notin X$}.
      \end{cases}
\]
Although we have not specified all the entries of the matrix, these entries show
that the matrix is triangular, and hence the determinant has absolute value
\[
\left(\Dfrac{\xi}{\pi-\xi}\right)^{|U|+|W|}
= e^{-O(\zeta |U|\log n)},
\]
by noting  $\xi\ge e^{-\zeta \log n} = o(1)$.
 \end{proof}

We are ready to bound the total contribution of regions $\varOmega_{U,W}$.
\begin{lemma}\label{outside}
$ \displaystyle
   \int_{\varOmega_3} |F(\thetavec)| d\thetavec = e^{-\omega(\log n )} J_0
$, where 
$\displaystyle \varOmega_3 :=
    \bigcup_{
	\substack{ 1\le |U|\le n/5 \\
	|W|\le \zeta^{-1} |U| } } \varOmega_{U,W}.$
\end{lemma}

\begin{proof}
 The volume of $\varOmega_2$ is only $e^{-\omega(n\log n)}$, so the bound
  $\abs{F(\thetavec)}\le 1$  and  Lemma~\ref{J0rough} give
\[ 
  \int_{\varOmega_2} \abs{F(\thetavec)}\,d\thetavec 
  = e^{-\omega(n\log n)}\, J_0.
\]
    Applying Lemma~\ref{JUWbound}  and bounding
   the number of choices of $W$ for given $U$  by $2^{|U|}$,  we obtain
\begin{align*}
   \sum_{ \substack{ 1\le |U|\le n/5 \\
	|W|\le \zeta^{-1} |U| } } 
   \int_{\varOmega_{U,W}-\varOmega_2} 
		   \abs{F(\thetavec)}\, d\thetavec
&= \sum_{ \substack{ 1\le |U|\le n/5 \\
	|W|\le \zeta^{-1} |U| } } e^{-\omega(|U| \log n)} 
		J_0 
		=  \sum_{ 1\le |U|\le n/5 } 2^{|U|} e^{-\omega(|U| \log n)}  J_0\\	 
&=        \sum_{t=1}^{n/5} \binom{n}{t} e^{-\omega(t\log n)}J_0
         =  \( \( 1 + e^{-\omega(\log n)}\)^n - 1\)J_0 \\
        &= O\( ne^{- \omega(\log n)}\) J_0 = O\(e^{- \omega(\log n)}\) J_0.
\end{align*}
This completes the proof of Lemma~\ref{outside}.
\end{proof}


\subsection{Proof of Theorem \ref{tm-eo}  and Corollary \ref{cor-eo}}\label{S:proof-main-main}

The value of $\EO(G)$ is given in~\eqref{NGint}. 
The main part of the integral is $J_0$, which is given to precision
$O(n^{-c})$ in Theorem~\ref{l:inside}. The remainder of the integral is
shown by Lemmas~\ref{S1bound}, \ref{l:UWcover} and \ref{outside}
to contribute the negligible additional error $O(e^{-\omega(\log n)})$.
This proves Theorem \ref{tm-eo}.

Corollary~\ref{cor-eo} follows from Theorem~\ref{tm-eo} on applying
the cumulant bounds in Theorem~\ref{T:cumulant_bounds} together
with Lemma~\ref{l:norms}(a).


\section{Regular tournaments and similar counts}\label{S:RT}

In this section we apply Theorem~\ref{tm-eo} to obtain terms
of the asymptotic expansion of the number of regular tournaments
$\RT(n)$.
We also use an entirely analogous theory to expand the
number $\ED(n)$ of Eulerian digraphs and the number
$\EOG(n)$ of Eulerian oriented graphs.
In all these cases, loops and parallel edges are forbidden,
as are 2-cycles except in the case of $\ED(n)$.
In the case of $\RT(n)$, $n$ must be odd.

The leading terms of these expansions were found by
McKay~\cite{mckay1990asymptotic}.

\begin{thm}\label{RT-ED-EOG}
For $n \to \infty$, with odd $n$ for $\RT(n)$,
\begin{align*}
\RT(n)
&= 
n^{1/2}
\Bigl(\Dfrac{2^{n+1}}{\pi n}\Bigr)^{(n-1)/2}
\exp\Bigl(
- \Dfrac{1}{2} + \Dfrac{1}{4n} + \Dfrac{1}{4n^2} 
+ \Dfrac{7}{24n^3} + \Dfrac{37}{120n^4} 
+ \Dfrac{31}{60n^5} + \Dfrac{81}{28n^6}  \\[0.6ex]
&{~}+ \Dfrac{5981}{336n^7} + \Dfrac{22937}{240n^8}
+ \Dfrac{90031}{180n^9} + \Dfrac{1825009}{660n^{10}} 
+ \Dfrac{4344847}{264n^{11}} + \O{ n^{-12} }
\Bigr), \\[0.6ex]
\ED(n)
&=
n^{1/2}
\Bigl(\Dfrac{4^{n}}{\pi n}\Bigr)^{(n-1)/2}
\exp\Bigl( - \Dfrac14 + \Dfrac{3}{16n} + \Dfrac{1}{8n^2} 
+ \Dfrac{47}{384n^3} + \Dfrac{371}{1920n^4} 
+ \Dfrac{1807}{3840n^5} + \Dfrac{655}{448n^6}  \\[0.6ex]
&{~}+ \Dfrac{435581}{86016n^7} + \Dfrac{1145941}{61440n^8} 
+ \Dfrac{13318871}{184320n^9} + \Dfrac{99074137}{337920n^{10}}
+ \Dfrac{1339710847}{1081344n^{11}} + \O{n^{-12}}
\Bigr), \\[0.6ex]
\EOG(n)
&= 
n^{1/2}
\Bigl(\Dfrac{3^{n+1}}{4\pi n}\Bigr)^{(n-1)/2}
\exp\Bigl( - \Dfrac{3}{8} + \Dfrac{11}{64n} + \Dfrac{7}{64n^2}
+ \Dfrac{233}{2048n^3} + \Dfrac{497}{2560n^4} + \Dfrac{27583}{61440n^5}  \\[0.6ex]
&{~}+ \Dfrac{55463}{43008n^6} 
+ \Dfrac{33678923}{7340032n^7} + \Dfrac{101414573}{5242880n^8}
+ \Dfrac{1882520759}{20971520n^9} + \Dfrac{101145677531}{230686720n^{10}}  \\[0.6ex]
&{~}+ \Dfrac{2469157786549}{1107296256n^{11}} +  \O{n^{-12}}
\Bigr).
\end{align*}
\end{thm}

As mentioned earlier, $\RT(n)$ is the special case of Theorem~\ref{tm-eo}
with $G=K_n$, corresponding to the constant term in 
$\prod_{j<k}\,(x_j/x_k+x_k/x_j)$.
Similarly, $\ED(n)$ is the constant term in
$\prod_{j< k} \,\((1 + x_j/x_k)(1+x_k/x_j)\)$ and $\EOG(n)$ is the
constant term in $\prod_{j<k}\,(1 + x_j/x_k+x_k/x_j)$.
This gives us integral representations as follows.

\begin{align*}
\RT(n)
&= \Dfrac{ 2^{n(n-1)/2} }{ (2\pi)^n }
\int_{U_{n}(\pi)}
\prod_{1 \le j < k \le n} \cos(\th_j - \th_k)\,d\thetavec \\
\ED(n) 
&= \Dfrac{ 2^{n(n-1)} }{ (2\pi)^n } 
\int_{U_{n}(\pi)}
\prod_{1 \le j < k \le n} 
\left( \dfrac12 + \dfrac12\cos(\th_j - \th_k) \right)\,d\thetavec \\
\EOG(n) 
&= \Dfrac{ 3^{n(n-1)/2} }{ (2\pi)^n }
\int_{U_{n}(\pi)}
\prod_{1 \le j < k \le n} 
\left( \dfrac13 + \dfrac23\cos(\th_j - \th_k) \right)\,d\thetavec.
\end{align*}

By Lemma~\ref{L:projections}(i), we are free to choose $w>0$
to form the covariance matrix $(L+wJ_n)^{-1}$ without changing
the cumulants, and since the off-diagonal terms are all the same
in these examples, we can choose $w$ so that $L+wJ_n$ is diagonal.
By symmetry, the diagonal entries will then all be equal, so we
are dealing with the cumulants of functions of a random vector
with i.i.d.\ Gaussian components.

For $\RT(n)$ we can choose $w=1$, so that the Gaussian vector
$\X$ has i.i.d.\ components
with distribution $\sqrt{\frac{n}{2\pi}}e^{- n x^2/2}$.
For $\ED(n)$ and $\EOG(n)$, we similarly obtain components
with distributions $\sqrt{\frac{n}{4\pi}}e^{- n x^2/4}$
and $\sqrt{\frac{n}{3\pi}}e^{- n x^2/3}$, respectively.

We continue the description for $\RT(n)$.
By Theorem~\ref{tm-eo}, when $c=12$, we can take $M=13$ and $K=14$.
In fact we can use $M=12$ since the logarithmic term in Lemma~\ref{l:norms}(a)
is unnecessary for $A=L+J_n$.
From \eqref{eq:rcumbnd2}, it is sufficient to choose $K=13$.

Since the connectivity condition in Theorem~\ref{cums} is
not convenient for exact computations, see \cite[Section 4.7.2]{zhang}, 
we found it more
efficient to compute the moments up to order 12.
The cumulants were then obtained by taking the logarithm of the moment generating function, equivalently to~\eqref{cumu},
\[
       \kappa_r(f_K(\X))
       = r!\, [t^r]
       \log\biggl(\,\sum_{i=0}^\infty
       \Dfrac{t^i}{i!}\,\E{ f_K(\X)^i}\biggr),
\]
where $[t^r]$ signifies extraction of the coefficient of $t^r$
in the formal expansion of the logarithm.
These moments are all $O(1)$, so it suffices to compute them
to additive error $O(n^{-12})$.

For $k \ge 0$, define the power sum symmetric functions
\[
\mu_k := \mu_k(\X) = \sum_{j=1}^n X_j^k.
\]
Then $\mu_0=n$. By noting that
\[
\E{ f_K(\X)^i }
= 
\E{ \left( \sum_{t\ge2}^{K}
c_{2t} \sum_{jk \in G}
(X_j -X_k)^{2t} \right)^{\!i\,}}
= 
\E{ \left( \dfrac12 \sum_{t\ge2}^{K}
c_{2t} \sum_{s=0}^{2t} \binom{2t}{s}
(-1)^s \mu_s \mu_{2t-s} \right)^{\!i\,}},
\]
the task then reduces to computing terms of the form
\[
\E {  \mu_{j_1}\cdots\mu_{j_m} },
\]
where $\{j_i\}_i$ might not be distinct.
This is the topic of the next section.


\subsection{Moments of i.i.d.\ Gaussian vectors}\label{S:bits}

Let $\X=(X_1,\ldots,X_n)$ be a Gaussian random vector with
independent and identically distributed components with each having mean 0 and variance~$n^{-1}$.
In this section, we provide some estimates of
the expectation of monomials in $\{\mu_k\}_k$.

\begin{lemma}\label{expbounds}
(a) For any integer $m\ge 0$,
\[
\E{ X_1^m } = \begin{cases}
0, & \text{ if $m$ is odd};\\
\dfrac{(m-1)!!}{ n^{m/2} }, & \text{ if $m$ is even.}
\end{cases}
\]
\\
(b) 
Suppose $m =O(1)$ and $j_i =O(1)$ for all $i\in[m]$.  Then
\[
 \E {  \mu_{j_1}\cdots\mu_{j_m} }
= \begin{cases}
0, & \text{ if $j_1+\cdots+j_m$ is odd};\\
\O{ n^{ J } }, & \text{ if $j_1+\cdots+j_m$ is even},
\end{cases}
\]
where
\[
J 
= J( j_1, \ldots, j_m )
= \Abs{\{k\st j_k\text{ is even}\}}
+ \dfrac{1}{2}\, \Abs{\{k\st j_k\text{ is odd}\}} 
- \dfrac{1}{2} \sum_{i \in [m]} j_i.
\]
\end{lemma}

\begin{proof}
Part (a) is the standard property of moments of Gaussian by noting $\sigma^2 (X_1) = 1/n$.
For (b),
first note that
\[
\E { \mu_{j_1}\cdots\mu_{j_m} } = 
\sum_{1 \le t_1,\ldots,t_m \le n}
\E{  X_{t_1}^{j_1} \cdots  X_{t_m}^{j_m} },
\]
in which $t_1,\ldots,t_m$ might not be distinct.
Since $X_1,\ldots,X_n$ are identical and independent,
then the expectation inside the sum depends only
on the equalities amongst the values $t_1,\ldots,t_m$.

Any particular sequence $t_1,\ldots,t_m$ defines a partition $\varPi(t_1,\ldots,t_m)$
with $q=q(t_1,\ldots,t_m)$ non-empty disjoint cells (sets), whose union is $\{1,\ldots,m\}$,
such that 
for any $1 \le a < b \le m$,
$t_a$ and $t_b$ are equal if and only if
they belong to the same cell. 

Then, by the independence of the $\{X_j\}$, we have
\begin{equation}\label{oneterm}
\E { X_{t_1}^{j_1} \cdots  X_{t_m}^{j_m} }
= \prod_{B\in\varPi(t_1,\ldots,t_m)} 
\E{ X_1^{\sum_{s\in B} j_s} },
\end{equation}
where the expectations on the right are provided by Lemma~\ref{expbounds}(a).

If $j_1+\cdots+j_m$ is odd, then 
there exists $B\in\varPi(t_1,\ldots,t_m)$ such that
$\sum_{s\in B} j_s$ is odd.
Hence from (a), we have that 
$\E { \mu_{j_1}\cdots\mu_{j_m} }  = 0$.

From now on, we assume that $j_1+\cdots+j_m$ is even.
For any $t_1,\ldots,t_m$ such that $\sum_{s\in B} j_s$ is even 
for all $B\in\varPi(t_1,\ldots,t_m)$,
we have that 
\begin{align}
\E { X_{t_1}^{j_1} \cdots  X_{t_m}^{j_m} }
&= \prod_{B\in\varPi(t_1,\ldots,t_m)} 
\E{ X_1^{\sum_{s\in B} j_s} } 
= \prod_{B\in\varPi(t_1,\ldots,t_m)} 
\left(\, \sum_{s\in B} j_s-1 \right)!!\, n^{- \frac{1}{2} \sum_{s\in B} j_s } \nonumber \\ 
&= n^{- \frac{1}{2} \sum_{i \in [m]} j_i} \prod_{B\in\varPi(t_1,\ldots,t_m)} 
\left(\, \sum_{s\in B} j_s-1 \right)!! \,
= \O{ n^{- \frac{1}{2} \sum_{i \in [m]} j_i} }.\label{evenSum}
\end{align}
It then suffices to count 
$t_1,\ldots,t_m$ such that $\sum_{s\in B} j_s$ is even 
for all $B\in\varPi(t_1,\ldots,t_m)$.
Let $( M_e, M_o )$ be a partition of $[m]$ 
such that $j_i$ is even for all $i \in M_e$
and $j_i$ is odd for all $i \in M_o$.
For each $i \in M_e$, we have $n$ choices of $t_i$.
For each $i \in M_o$, 
there must exist $j \in M_o$ such that
$i \ne j$ and $t_i = t_j$,
since 
$\sum_{s\in B} j_s$ is even 
for all $B \in\varPi(t_1,\ldots,t_m)$.
So we have $|M_o|/2$ pairs with each having $n$ choices.
Therefore, in total there are $|M_e| + |M_o|/2$ ways of choosing with each having $n$ options, 
then combining with \eqref{evenSum} completes the proof.
\end{proof}

Equation~\eqref{oneterm}, together with the fact that $\varPi(t'_1,\ldots,t'_m)
=\varPi(t_1,\ldots,t_m)$ for exactly $(n)_q$ index sets $t'_1,\ldots,t'_m$,
allows us to compute $ \E { \mu_{j_1}\cdots\mu_{j_m} }$
by summing over all partitions of $\{1,\ldots,m\}$.
However, this is very inefficient since the number of partitions 
(the Bell numbers) grow very quickly.

Note that many partitions give the same contribution.
We define $\Naturals=\{0,1,\ldots\,\}$ and $\Naturals_+=\{1,2,\ldots\,\}$.
We will represent a multiset as a pair $(S, \nu)$, 
where $S$ is a set that contains distinct elements 
and~$\nu : S\to\Naturals$ is a function
that specifies the number of times an element of~$S$ belongs to the multiset.
We will also use the notation $\ms{z_1,\ldots,z_k}$
for a multiset by listing all its elements.
For example, $\ms{2,2,2,3,3,4}$ is the multiset $(\Naturals_+,\nu)$
where $\nu(2)=3$, $\nu(3)=2$, $\nu(4)=1$, and $\nu(k)=0$ otherwise.

Now we introduce a definition that 
is useful for calculation.
\begin{df}
A \textit{cell type} is a multiset $(\Naturals_+,\nu)$.
A \textit{partition type} $(C,\eta)$ is a multiset of cell types, where
$C$ is the set of all cell types.
\end{df}

Then
the expansion of $\mu_{j_1}\cdots\mu_{j_m}$ contains $n^m$ terms
such that each term $X_{t_1}^{j_1}\cdots X_{t_m}^{j_m}$
has a unique partition type where each cell type corresponds to the exponents
of a set of equal indices.
For example, in the expansion of $\mu_1\mu_1\mu_2\mu_2\mu_4$, the term
$X_5^1X_5^1 X_6^2 X_5^2 X_6^4$ has
partition type $\Ms{\ms{1,1,2},\ms{2,4}}$, since $X_5^1,X_5^1,X_5^2$ have the
same index and $X_6^2, X_6^4$ have the same index.

Note that the multiset union of the cell types in the partition type equals the multiset
of the indices of the monomial $\mu_1\mu_1\mu_2\mu_2\mu_4$, namely,
$\ms{1,1,2,2,4}$.
Also note that for instance, 
the different term $X_3^1X_3^1 X_3^2 X_2^2 X_2^4$ in
the expansion of $\mu_1\mu_1\mu_2\mu_2\mu_4$ has the same
partition type even though the positions of the equal indices changed
from $\{1,2,4\},\{3,5\}$ to $\{1,2,3\},\{4,5\}$.

For $k\in\Naturals_+$, let $\widehat\nu(k)$ be the number of times~$\mu_k$
appears in $\mu_{j_1}\cdots\mu_{j_m}$.
Clearly there is no term in the expansion of  $\mu_{j_1}\cdots\mu_{j_m}$ with
partition type~$(C,\eta)$ unless
\begin{equation}\label{ptncond}
\sum_{(\Naturals_+,\nu)\in (C,\eta)} \nu(k) =\widehat\nu(k),
\end{equation}
for all $k\in\Naturals_+$.

\begin{lemma}
We have that
\[
\E { \mu_{j_1}\cdots\mu_{j_m} }
= \sum_{(C,\eta)} 
A_{(C,\eta)} B_{(C,\eta)}  \prod_{\tau\in C} \,\E{  X_1^{\sum_{j\in \tau} j }},
\]
where the sum is over all partition types satisfying~\eqref{ptncond},
\[
A_{(C,\eta)}
:= 
\Biggl( \biggl( n-\sum_{\tau\in C} \eta(\tau) \biggr)! \prod_{\tau\in C} \eta(\tau)! \Biggr)^{\!-1}
n!,
\]
and
\[
B_{(C,\eta)}
:= \prod_{k\in\Naturals_+} 
\Biggl(\, \prod_{(\Naturals_+,\nu)\in (C,\eta)} \nu(k)! \Biggr)^{\!-1} 
\Biggl(\, \sum_{(\Naturals_+,\nu)\in (C,\eta)} \nu(k) \Biggr)!.
\]
\end{lemma}
\begin{proof}
Suppose $(C,\eta)$ is a partition type satisfying~\eqref{ptncond}.
Since some of the cells in ${(C,\eta)}$ are the same, 
then the number of ways to assign distinct indices from $\{1,\ldots,n\}$ is $A_{(C,\eta)}$.
The number of ways to assign the $m$ positions of $\mu_{j_1}\cdots\mu_{j_m}$
to the cells of~${(C,\eta)}$ is $B_{(C,\eta)}$. 
This completes the proof.
\end{proof}

Here we include an illustrative example.
If ${(C,\eta)} = \Ms{ \ms{1, 1, 2}, \ms{1, 1, 2}, \ms{1, 2, 2} }$,
then 
$C = \bigl\{ \ms{1, 1, 2}, \ms{1, 2, 2} \bigr\}$
and we have
$\eta( \ms{1, 1, 2} ) = 2$ and $\eta( \ms{1, 2, 2} ) = 1$.
So 
\[
A_{(C,\eta)} = \Dfrac{n!}{(n-3)!\,2!\,1!} = \dfrac{1}{2} n(n-1)(n-2),
\]
since we need to choose three distinct indices from $[n]$ for each cell type 
such that two of them are of the same type.
We also have 
\[
B_{(C,\eta)} = \dfrac{ 5! }{ 2!\, 2! } \dfrac{4!}{2!},
\]
since 
we need to choose the exponents for each cell type,
and that is $\binom{5}{2}\binom{3}{2}$ for one
and $\binom{4}{2}\binom{2}{1}$ for two.

This is the key to our method because the contribution of a term
to $\E { \mu_1\mu_1\mu_2\mu_2\mu_4 } $ only depends on its
partition type~$(C,\eta)$, namely
\[
\prod_{\tau\in C} \,\E{ X_1^{\sum_{j\in \tau} j} }.
\]
Thus, to evaluate $\E{ \mu_{j_1}\cdots\mu_{j_m} }$, we can sum over
partition types, rather than over partitions.
This is a very large improvement, since the number of partition types
is much smaller.
As an example, for $\mu_2^{10}\mu_3^2\mu_4^{10}$ the number of partitions
is 4,506,715,738,447,323, but the number of partition types is only 360,847.

In Table~\ref{tab:moments} we give the leading terms of
$\E{f_K(\X)^r}$ and $\kappa_r(f_K(\X))$ for regular tournaments
for small~$r$.

\begin{table}[h!]
\centering
\renewcommand{\arraystretch}{1.5} 
\begin{tabular}{c|c|c}
\hline
$r$ & $\E{ f_5(\mathbf{X})^r}$ & $\kappa_r(f_5(\mathbf{X}))$ \\ \hline
1 & $-\frac{1}{2} - \frac{5}{6n} - \frac{13}{3n^2} - \frac{137}{5n^3} + O(n^{-4})$
& $-\frac{1}{2} - \frac{5}{6n} - \frac{13}{3n^2} - \frac{137}{5n^3} + O(n^{-4})$\\ 
2 &  
$\frac{1}{4} + \frac{3}{n}
+ \frac{811}{36n^2}
+ \frac{8429}{45n^3} + O(n^{-4})$
  & $\frac{13}{6n} + \frac{35}{2n^2} + \frac{6871}{45n^3} + O(n^{-4})$ \\ 
3 &  
$-\frac{1}{8} - \frac{31}{8n}
- \frac{1463}{24n^2}
- \frac{162701}{216n^3} + O(n^{-4})$
  & $-\frac{25}{n^2} - \frac{1261}{3n^3} + O(n^{-4})$ \\ 
4 &  
$\frac{1}{16} + \frac{11}{3n}
+ \frac{835}{8n^2}
+ \frac{559643}{270n^3} + O(n^{-4})$
  & $\frac{1541}{3n^3} + O(n^{-4})$ \\
5 & 
$-\frac{1}{32} - \frac{95}{32n}
- \frac{9745}{72n^2}
- \frac{890933}{216n^3} + O(n^{-4})$
  & $O(n^{-4})$ \\\hline
\end{tabular}
\caption{
$\E{ f_5(\mathbf{X})^r}$ 
and $\kappa_r(f_5(\mathbf{X}))$
for regular tournaments\label{tab:RT}}
\label{tab:moments}
\end{table}


\nicebreak
\section{Tail bounds for cumulant expansions}
\label{S:cumu-tm}

We prove Theorem \ref{cumuThm} in this section.
Let $\X=(X_1,\ldots, X_n)$ be a random vector with independent components
taking values in some product space $\S := S_1\times\cdots\times S_n$.
Let $\calF_{\infty}(\X)$ be the space of bounded real functions on $\S$,  measurable with respect to $\X$, equipped with the infinity norm
\[
\norm{f}_{\infty}  = \sup_{\xvec \in \S} |f(\xvec)|, \qquad  f \in \calF_{\infty}(\X).
\]
For a linear operator $F$ on $\calF_{\infty}(\X)$,  we consider the standard induced operator norm
\[
\norm{F}_{\infty} :=  \sup_{ \substack{f\in\calF_{\infty}(\X) \\ \norm{f}_\infty>0}} \frac{\norm{F[f]}_\infty}{\norm{f}_\infty}.
\]

\subsection{Expectation and difference operators}

For any $V \subseteq [n]$,  define the operator $\Ex^{V}$ on $\calF_\infty(\X)$ by 
\[
	\Ex^{V}[f](\xvec):= \Ex\left[f(\X) \mid X_j = x_j \text{ for $j \notin V$}\right],
\]
 Informally,  $\Ex^{V}$ corresponds to ``averaging'' with respect to all $X_j$ with $j \in V$.
  Since  expectation can not exceed the supremum, we obtain 
  $
  	\norm{\Ex^{V}}_\infty \leq 1. 
  $
  Observe also that  the operators $\Ex^{V}$ and $\Ex^{V'}$ commute for any $V, V'\subseteq [n]$:
  \[
  	\Ex^{V} \Ex^{V'} = 	\Ex^{V'} \Ex^{V} = \Ex^{V\cup V'}.
  \]
If $V = \{j\}$, we write $\Ex^{j}:=\Ex^{\{j\}}.$

Recall from \eqref{def:R-delta} that 
\[
 R^{V}_{\yvec}:= R^{v_1}_{\yvec}\cdots R_{\yvec}^{v_k}, \qquad \partial^V_{\yvec}:= \partial_{\yvec}^{v_1}\cdots \partial_{\yvec}^{v_k},
\]
where  $\partial_{\yvec}^j:= I - R^j_{\yvec}$ and  $I$ is the identity operator.  This definition does not depend on the order of elements in $V$ since the operators $ R^j_{\yvec} $  and $ R^{j'}_{\yvec} $ commute for any $j,j'\in [n]$. 
Clearly $\norm{R_{\yvec}^V}_\infty \leq 1$. We also have that
\begin{equation}\label{RW-commute}
\partial_{\yvec}^{V} R_{\yvec}^{V'} = 	R_{\yvec}^{V'} \partial_{\yvec}^{V}.
\end{equation}
We set  $R^{\emptyset}_{\yvec} 
= \partial^{\emptyset}_{\yvec} = I$.
We will use the identity
\begin{equation}\label{delta-identity}
 \partial_{\yvec}^V  = \sum_{W\subset V} (-1)^{|W|} R_{\yvec}^W.
\end{equation}
From \eqref{def:Delta}, we have that  
\begin{equation*} 
	 \Delta_{V}(f)   =  \sup_{ \y\in \S} \|\partial^V_{\yvec}[f]\|_\infty.
\end{equation*}
In particular,
$\Delta_{\emptyset}(f) := \sup_{ \xvec \in \S} \left|f (\xvec)\right|.$
By definition, $\Delta_V$ satisfies the triangle inequality
\begin{equation}\label{Delta:triangle}
		\Delta_V(f + g) \leq \Delta_V(f)+\Delta_V(g).
	\end{equation}

\begin{lemma}\label{L:E-Delta}
	 For any $f \in \calF_{\infty}(\X)$,    $V \subseteq [n]$,  and $j \in [n] $, we have
		\[
		 \Delta_{V} (\Ex^j[f]) \leq \Delta_V(f), \qquad  
		 \Delta_{V} (f - \Ex^j[f]) \leq \Delta_{V\cup \{j\}}(f).
			\]
			Furthermore, if $j \in V$ then $\Delta_V(\Ex^j [f]) = 0$ and $\Delta_{V} (f - \Ex^j[f]) =  \Delta_{V}(f)$.
	\end{lemma}
\begin{proof} 
	We start with the second part, which is the case when $j\in V$.  For any $\yvec \in \S$, we have that 
$
		R_{\y}^j \Ex^j [f] = \Ex^j [f]
$
	because $\Ex^j [f] (\xvec)$ does not depend on $j$th component of $\x$. This implies that
	\[
		\partial_{\y}^V \Ex^j [f] = 0 \quad \text{ and } \quad \partial_{\y}^V [f - \Ex^j [f]] = \partial_{\y}^V [f].
	\]
	The second part follows.
	
	Now we assume that $j\notin V$.
	Let $\yvec \in \S$ 
 and $\Y = (y_1,\ldots,y_{j-1}, X_j, y_{j+1}, y_n)$. 
 For any $W\subseteq[n]\setminus \{j\}$, we have
	\[
		R_{\yvec}^W[f] = \Ex\left[R_{\Y}^W [f] \right] \quad \text{ and }\quad (R_{\yvec}^W \Ex^j) [f] = \Ex \left[ R_{\Y}^{W\cup\{j\}}[f]\right].
	\]
	Then, using \eqref{delta-identity}, we find that
	\begin{equation*}
			\partial_{\yvec}^V [\Ex^j [f]] =  \sum_{W\subseteq V} (-1)^{|W|} (R_{\yvec}^W \Ex^j)[f]
			= \Ex \left[ \sum_{W\subseteq V} (-1)^{|W|}  R_{\Y}^{W\cup\{j\}}[f]] \right] = \Ex\left[\partial_{\Y}^W R^j_{\Y}[f]\right]
   \end{equation*}
and 
\begin{align*}
		\partial_{\yvec}^V [f - \Ex^j [f]] &= \sum_{W\subseteq V} (-1)^{|W|} R_{\yvec}^W[f- \Ex^j [f]]
		\\&= \Ex \left[ \sum_{W\subseteq V} (-1)^{|W|}  \left(R_{\Y}^W[f] - R_{\Y}^{W\cup\{j\}}[f]] \right) \right] 
		= \Ex \left[ \partial_{\Y}^{V\cup\{j\}} [f]\right].
	\end{align*}
 Since the expectation can not exceed the supremum, we obtain 
 \begin{align*}
 			\partial_{\yvec}^V [\Ex^j [f]] &\leq \sup_{\y\in\S}\|\partial_{\yvec}^{V}[f]\|_\infty = \Delta_{V}(f),\\
     	\|\partial_{\yvec}^V [f - \Ex^j [f]]\|   &\leq \sup_{\y\in\S}\|\partial_{\yvec}^{V\cup\{j\}}[f]\|_\infty = \Delta_{V\cup\{j\}}(f).
 	\end{align*}
  Taking the supremum over $\yvec$ completes the proof.
	 	\end{proof}

Let $D_r(V)$ denote the set of all dissections of $V\subseteq [n]$ into an ordered collection of $r$  subsets 
$(V_1,\ldots, V_r)$,  that is, the sets $V_j$ are disjoint (possibly empty) and $V = V_1 \cup \cdots \cup V_r$.
 \begin{lemma}\label{lemma:Deltabound}
 	Let $f_1, \ldots, f_r \in \calF_{\infty} (\X)$ and $V \subseteq [n]$.
 	Then 
 	\[
 	 \Delta_V (f_1\cdots f_r) \leq 
 	 \sum_{ (V_1,\ldots, V_r)\in D_r(V)}\,
 	\prod_{j=1}^r\,  \Delta_{V_j}(f_j).
 	\]
 \end{lemma}

 \begin{proof}
 	The statement is trivial for $r=1$. We proceed to the case when $r=2$.
     	For any $\yvec\in \S$,  $j  \in [n]$,	   observe that
 	\begin{equation*} 
 		 \partial^{j}_{\yvec}[f_1 f_2] = f_1 \cdot  \partial_{\yvec}^j[f_2] +   \partial_{\yvec}^j[f_1] \cdot R_{\yvec}^{j}[f_2].
 	\end{equation*}
 Applying this analog of the product rule of differentiation several times and using \eqref{RW-commute}, we find that
 \[
 	\partial^V_{\yvec}[f_1 f_2] = \sum_{W \subseteq V} \partial_{\yvec}^W[f_1] \cdot R_{\yvec}^W[ \partial_{\yvec}^{V\setminus W}f_2].
 \]
 Recalling that $\|R_{\yvec}^W\| \leq 1$ and using  the triangle inequality and definition \eqref{def:Delta}, we find  that
 \[
 	\|\partial^V_{\yvec}[f_1 f_2]\|_\infty \leq \sum_{W \subseteq V} \| \partial_{\yvec}^W[f_1]\|_\infty \cdot \| \partial_{\yvec}^{V\setminus W}[f_2]\|_\infty
 	 \leq \sum_{W \subseteq V} \Delta_{W}(f_1) \Delta_{V\setminus W}(f_2).
 \]
Taking the supremum over $\yvec$ completes the proof for the case when $r=2$.  

The statement for $r>2$ follows from the bound above for $r=2$ by bounding 
\[
		 \Delta_V (f_1\cdots f_r) \leq  \sum_{W \subseteq V}  \Delta_W(f_1) \Delta_{V\setminus W} (f_2\cdots f_r)
\]
and using a simple inductive argument.
 \end{proof}

\nicebreak
 \subsection{Cumulant identities and bounds} 
 For each $j\in [n]$, 
 let 
\[
\Ex^{\geq j}  :=  \Ex^{(\{j,\ldots,n\})} .
\]
 For  $f_1,\ldots,f_t\in \calF_\infty(\X)$,  consider the conditional joint cumulant  defined by
 \begin{equation}\label{def:cum-operator}
 	\kappa^{\geq j}[f_1, \ldots, f_r] = \sum_{\tau \in \mathcal P_r} \, (\card{\tau}-1)!\, (-1)^{\card{\tau}-1}
 	\prod_{B\in\tau} \Ex^{\geq j}\Bigl[\prod_{k \in B} f_k\Bigr],
 \end{equation}
 where  $\mathcal P_r$ denotes the set of   unordered partitions $\tau$ of~$[r]$ (with non-empty blocks)  and $|\tau|$ denotes  the number of blocks in the partition $\tau$. 
 We also set
 \[ \Ex^{(\geq n+1)}[f]   = \kappa^{(\geq n+1)}[f]:= f \ \ \  \text{ and } \ \
\kappa^{(\geq n+1)}[f_1,\ldots,f_r] :=  0 
\text{ for } r \geq 2.
\]

 \begin{lemma}\label{cumulantgeneral}
Let $f_1, \ldots, f_r \in \calF_{\infty}(\X)$. The following hold.
 	\begin{itemize} \itemsep=0pt
 		\item[(a)]   $	\kappa^{\geq j}$ is a symmetric function and also a multilinear function, that is, 
 		\[
 			\kappa^{\geq j}[a_1 f_1 + a_2 f_1', f_2, \ldots, f_r] = a_1 	\kappa^{\geq j}[f_1,\ldots, f_r] + a_2 	\kappa^{\geq j}[ f_1', f_2,\ldots, f_r]
 		\]
 		for any $a_1,a_2 \in \Reals$ and $f_1'  \in \calF_{\infty}(\X)$. Furthermore, if
 		 $r\geq 2$ and $f_r = \Ex^{\geq j} g$ for some  $g \in \calF_{\infty}(\X)$,  then 
 		$	\kappa^{\geq j}[f_1,\ldots,f_r]\equiv 0$.
 		
 		\item[(b)] Let $\log(1+t) = \sum_{k=1}^\infty \frac{(-1)^k}{k} t^k $. We have
 	    \[ 	\kappa^{\geq j}[f_1,\ldots, f_r]  = [t_1\cdots t_r]\log \biggl(1+\sum_{k=1}^\infty 
 		\frac{	\kappa^{\geq j}[t_1f_1+\cdots+t_r f_r]^k}{k!} \biggr),\]
 		where $t_1,\ldots,t_r$ are real indeterminants, and $[t_1\cdots t_r]$ indicates
 		coefficient extraction in the formal series expansion. 
 		\item[(c)]  Let $k \in \{j,\ldots ,n+1\}$.  Then
 		\[
 		\kappa^{\geq j}[f_1,\ldots,f_r] = 
 		\sum_{p=1}^r 
   \sum_{\{B_1,\ldots,B_p\}\in \mathcal P_r} 
 		\kappa^{\geq j}\left[\kappa^{\geq k}[f_j:j\in B_1],\ldots,
 		\kappa^{\geq k}[f_j:j\in B_p]\right],
 		\]
 		where  
 		\[
 		\kappa^{\geq k}[f_j:j\in \{i_1,\ldots, i_{\ell}\} ] := \kappa^{\geq k}[f_{i_1}, \ldots, f_{i_\ell}].
 		\]
 		\item[(d)] 
 		For any set $V \subseteq [j-1]$, we have
 		\[
 		 \Delta_{V}  \(\kappa^{\geq j}[f_1, \ldots, f_r]\) 
 		\leq  \(\dfrac32\)^r   (r{-}1)! \!\sum_{(V_1,\ldots, V_r) \in D_r(V)} \,
 		\prod_{k=1}^r \, \Delta_{V_k} (f_k ).
 		\]
 	\end{itemize}
 \end{lemma}
 \begin{proof}
 	The fact that $\kappa^{\geq j}$ is symmetric and multilinear follows immediately
 	from the definition.
 	We proceed to the second part of (a).  Consider the terms of the defining summation in \eqref{def:cum-operator} that correspond to the partition $\tau'$ of $[r{-}1]$ that results
 	from disregarding~$r$.  
 	For any partition $\tau'=\{B_1,\ldots,B_k\}$, there are exactly $k+1$ corresponding terms.
 	One has~$r$ by itself, $\tau_0=\{\{r\},B_1,\ldots,B_k\}$ with
 	coefficient $(-1)^k k!$, and $k$ have the form
 	$\tau_j=\{B_1,\ldots,B_{j-1},B_j\cup\{r\},B_{j+1},\ldots,B_k\}$ with
 	coefficient $(-1)^{k-1}(k-1)!$.
	Moreover for $0\le j\le k$ we have
\[ 
\prod_{B\in \tau_j} \Ex^{\geq j}\Bigl[ \prod_{k\in B} f_k \Bigr]
 	= f_r \prod_{i=1}^k \Ex^{\geq j}\Bigl[ \prod_{k\in B_i} f_k \Bigr].
\]
	Since the coefficients have zero sum,
 we have $\kappa^{\geq j}[f_1,\ldots,f_r]\equiv 0$.
 	
 	Parts (b) and (c) are proved by Speed~\cite{MR725217} for random
 	variables when $F$ is
 	the expectation operator and $G$ is a conditional expectation
 	operator.  It is easy to check that, in addition to the combinatorial
 	properties of the partition lattice, only the linearity of expectation
 	and the law of total expectation are used,
 	so the same proofs work here also.

For (d), combining triangle inequality \eqref{Delta:triangle} and  definition   \eqref{def:cum-operator}, we obtain
 	\[
 	 \Delta_V(\kappa^{\geq j}[f_1, \ldots, f_r]) 
 	  \leq \sum_{\tau \in \mathcal P_r} \, (\card{\tau}-1)! \,
 	  \Delta_V \left(\prod_{B\in\tau} \Ex^{\geq j}\Bigl(\prod_{k\in B} f_k\Bigr)\right).
 	\]
 Applying Lemma \ref{L:E-Delta} several times, we obtain,  for any $W\subseteq V$,
 \[
 	 \Delta_{W} \left(\Ex^{\geq j}\left[\, \prod_{j\in B} 
 	f_j \right]\right) \leq  \Delta_{W} (f_j).
 \]
Then,  using  Lemma \ref{lemma:Deltabound} twice, we obtain
 \begin{align*}
 	\Delta_V \left(\prod_{B\in\tau} \Ex^{\geq j}\left[\,\prod_{j\in B} f_j\right]\right) &\leq
 	\sum_{(U_B)_{B\in \tau} \in D_{|\tau|}(V)} \prod_{B\in \tau} \Delta_{U_B} \left(\Ex^{\geq j}\left[\, \prod_{k\in B} 
 	  f_k \right]\right)\\
 	  &\leq 
 	  \sum_{(U_B)_{B\in \tau}\in D_{|\tau|}(V)} \prod_{B\in \tau} \Delta_{U_B} \left(\, \prod_{k\in B} 
 	  f_k \right)
 	  \leq \sum_{(V_1,\ldots, V_r) \in D_r(V)}   
 	  \prod_{k=1}^r  \Delta_{V_k} (f_k ).
 \end{align*}
From Lemma~\ref{ordpart},
\[
 \sum_{\tau\in \mathcal P_r} (\abs\tau-1)!
 \le r! \sum_{k=1}^r \Dfrac{2^{k-r}}{k}\binom{r-1}{k-1}
 =(r-1)!\,2^{-r}(3^r-1) \le \bigl(\dfrac32\bigr)^r(r-1)!,
\]
completing the proof.
 \end{proof}

The conditional cumulant of order $r$ is defined by
 \[
 \kappa_r^{\geq j}[f] =  \kappa^{(\geq j)}[\underbrace{f,\ldots, f}_{r \text{ times}}].
 \]
 Applying Lemma \ref{cumulantgeneral}, we derive the following properties of $\kappa_r^{\geq j}$.

 \begin{lemma}\label{l:cum_rec} 
 	If $f \in \calF_{\infty}(\X)$, $j \in [n]$,  and integer $r \ge 1$, then  the following hold. 
 	\begin{itemize}
 		\item[(a)]  For any $k\in[n+1]$, $k\geq j$, we have
 		\[
 		\kappa^{\geq j}_r [f]
 		=  \sum_{p=1}^r  \sum_{ \{B_1,\ldots,B_p\}\in \mathcal P_r} 
 		\kappa^{\geq j} \left[ \kappa^{\geq k}_{|B_1|}[f],
 		\ldots,
 		\kappa^{\geq k}_{|B_p|}[f]\right].
 		\]   
 		
 		\item [(b)] For any $V \subseteq [j-1]$ and $r\geq 2$, we have
 		\[
 			\Delta_V  \left(\kappa^{\geq j}_r [f]\right) \leq  \sum_{k=j}^n\, \sum_{p=2}\limits^r
 			\sum_{\{B_1,\ldots,B_p\}\in \mathcal P_r}   \(\dfrac32\)^p   (p{-}1)!  
			\sum_{(V_1,\ldots,V_p) \in D_p(V)}      
 			\prod_{r=1}^{p} \Delta_{V_r \cup\{k\}} \left(\kappa^{\geq k+1}_{|B_r|} [f]\right),
 		\]

 	\end{itemize}  
 	
 \end{lemma}
 \begin{proof}
 	Part  (a) is just a special case of 	Lemma \ref{cumulantgeneral}(c). For (b),  applying part (a) with  $k= j+1$  and using triangle inequality \eqref{Delta:triangle},  we obtain
 	 	\[
 		\Delta_V \left(\kappa^{\geq j}_r [f]\right)
 		 \le  
		 \Delta_{V} \left(\Ex^{\geq j} \kappa^{\geq j+1}_r [f]\right)
+ \sum_{p=2}^r
 		\sum_{ \{B_1,\ldots,B_p\}\in \mathcal P_r} 
 		\Delta_V \left(\kappa^{\geq j} \left[ \kappa^{\geq j+1}_{|B_1|}[f],
 		\ldots,
 		\kappa^{\geq j+1}_{|B_p|}[f]\right]\right).
 	\]
 From  Lemma \ref{L:E-Delta}, we know that
 	\[
 	  \Delta_V \left(\Ex^{\geq j} \kappa^{\geq j+1}_r [f] \right)\leq 
 	  \Delta_V \left( \kappa^{\geq j+1}_r [f] \right).
 	\]
 	Using Lemma \ref{cumulantgeneral}(a), we find that 
 	\begin{equation}\label{cumulant-red}
 	\begin{aligned}
 		\kappa^{\geq j} \left[ \kappa^{\geq j+1}_{|B_1|}[f],
 		\ldots,
 		\kappa^{\geq j+1}_{|B_p|}[f]\right]
 		&=                  
 	\kappa^{\geq j} \left[ (I - \Ex^j)\kappa^{\geq j+1}_{|B_1|}[f],
 	\kappa^{\geq j+1}_{|B_2|}[f],
 	\ldots,
 	\kappa^{\geq j+1}_{|B_p|}[f]\right]
 		\\
 		\cdots &=               	\kappa^{\geq j} \left[ (I - \Ex^j)\kappa^{\geq j+1}_{|B_1|}[f],
 		\ldots,
 		(I -\Ex^j)\kappa^{\geq j+1}_{|B_p|}[f]\right].
 	\end{aligned}
 \end{equation}
 	Then, applying Lemma \ref{cumulantgeneral}(d) and using Lemma \ref{L:E-Delta}
 	to estimate 
 	\[ 
 	\Delta_{V_r} \left( (I-\Ex^j)\kappa^{\geq j+1}_{|B_r|} [f] \right) \leq 
 	\Delta_{V_r \cup\{j\}} \left(\kappa^{\geq j+1}_{|B_r|} [f] \right),
 	\]
  we obtain  
 	\begin{align*}
 	\Delta_V \left(\kappa^{\geq j}_r [f] \right)&\leq  \Delta_V \left(\kappa^{\geq j+1}_r [f] \right)
 		\\
 		&\qquad+
 		\sum_{p=2}\limits^r
 		\sum_{\{B_1,\ldots,B_p\} \in \mathcal P_r} \(\dfrac32\)^p   (p{-}1)!  \sum_{(V_1,\ldots,V_p) \in D_p(V)}      
 		\prod_{r=1}^{p} \Delta_{V_r \cup\{j\}} \left(\kappa^{\geq j+1}_{|B_r|} [f] \right).
 	\end{align*}
 	Estimating similarly  $	\Delta_V \left(\kappa^{\geq j}_r [f] \right)$ for $k=j+1,\ldots, n$
 	and recalling that $\kappa_r^{\geq n+1} [f] \equiv 0$ for $r\geq 2$,  we have proved part (b).
 \end{proof}
 
 
 \subsection{Estimates when  the sums  of $\Delta_V$  are not large}

 For $v \in [n]$, let  
 \[
S_v (f):= \max_{j \in [n]}  \sum_{V\in \binom{[n]}{v} \st j \in V } \Delta_V(f).
 \]

 Throughout this section, we assume that $S_v (f)$ is not very big.  Namely,  for $\alpha \geq 0$ and positive integer $m$, let 
 \[
 		\calF^{\alpha}_m(\X):= \left\{f \in \calF_\infty(\X) \st S_v (f)  \leq \alpha    \text{ for all $v \in [m]$}\right\}.
\]
 In particular, for any $f \in \calF^{\alpha}_m(\X)$, we have
 \[
 	\max_{j\in [n]} \Delta_j(f) =  S_1(f) \leq \alpha.
 \]
 
\begin{lemma}\label{l:d-cumulant1}
	Suppose $f \in 	\calF^{\alpha}_m(\X)$ for some   $\alpha \geq 0$ and positive integer $m$.
	Then, for any $r \in [m]$ and  $j\in [n]$,
	we have 
	$$
	   \left\|\kappa_r^{\geq {j+1}}[f] - \Ex^j \kappa_r^{\geq {j+1}}[f]\right\|_\infty \leq   80^{r-1} \frac{(r-1)!}{r} \alpha^r.
	$$
\end{lemma}
\begin{proof}
	First, recalling $\kappa_1^{\geq j+1} = \Ex^{\geq j+1}$ and using   Lemma \ref{L:E-Delta},
   we find for any $V \subseteq [j]$ that
	\begin{equation}\label{eq:E-Delta}
 \Delta_{V} \left( \kappa_1^{\geq j+1} [f]\right)
 =\Delta_{V} \left( \Ex^{j+1} \cdots \Ex^n [f]\right)
  \leq  \Delta_V(f). 
\end{equation}
From Lemma \ref{L:E-Delta}, we also find that if $V \cap \{j+1 \ldots n\} \neq \emptyset$ then 
\begin{equation}\label{eq:kappa-zero}
	\Delta_V \left(\kappa_r ^{\geq j+1}[f]\right) = \Delta_V \left( \Ex^{\geq j+1}\kappa_r ^{\geq j+1}[f]\right) =0.
\end{equation}

	We prove the following statement by induction on $r \in [m]$: for any $j\in  [n]$  and $ v \in [m-r+1]$, we have
	\begin{equation}\label{ind:hyp}
		S_v(\kappa_{r}^{\geq j+1}[f]) \leq \hbar(v, r):= 56^{r-1} \frac{(r-1)!}{r} \binom{v+2r-3}{r-1} \alpha^r.
	\end{equation}
	 If $r=1$ then we see from \eqref{eq:E-Delta} 
	and \eqref{eq:kappa-zero} that 
	\[
		S_v(\kappa_1^{\geq j+1}[f]) \leq S_v[f]   \leq \alpha = \hbar(v, 1).
	\]
Thus, we verified the base of induction.

For the induction step,    from  Lemma \ref{l:cum_rec}(b), for any $V \in \binom{[n]}{v}$, we find that
\[
 \Delta_V \left(\kappa_r^{\geq j+1} [f]  \right) 
	\leq \sum_{k=j+1}^n\,
	\sum_{p=2}^r
	\sum_{
		\{B_1,\ldots,B_p\}\in \mathcal P_r
	}   \(\dfrac32\)^p   (p{-}1)! \hspace{-4mm}  \sum_{(V_1,\ldots,V_p)\in D_p(V)}      
	\prod_{t=1}^{p} \Delta_{V_p \cup\{k\}} \left(\kappa^{\geq k+1}_{|B_t|} [f] \right).
\]
Applying the induction hypothesis,   we find that, for any $ i   \in[j] $,
\begin{align*}
	\sum_{V \in \binom{[n]}{v} : i\in V} &   \sum_{\substack{(V_1,\ldots,V_p)\in D_p(V)\\ i \in V_1}}   \sum_{k=j+1}^n   
	\prod_{t=1}^{p} \Delta_{V_p \cup\{k\}} \left(\kappa^{\geq k+1}_{|B_t|} [f] \right) 
	\\
		&\leq  \sum_{\substack {v_1,\ldots, v_p \in \mathbb{N} \\ v_1 +\cdots +v_p = v}} 
		\sum_{k=j+1}^n  \sum_{\substack{V_1 \in \binom{[n]}{v_1+1} \\ i, k \in V_1} }
		\sum_{\substack{V_2 \in \binom{[n]}{v_2+1} \\  k \in V_2} } \cdots 
		\sum_{\substack{V_p \in \binom{[n]}{v_p+1} \\  k \in V_p} }
		\prod_{t=1}^{p} \Delta_{V_p \cup\{k\}} \left(\kappa^{\geq k+1}_{|B_t|} [f] \right)
	\\
	&\leq  \sum_{\substack {v_1,\ldots, v_p \in \mathbb{N} \\ v_1 +\cdots +v_p = v}} 
	\sum_{k=j+1}^n  \sum_{\substack{V_1 \in \binom{[n]}{v_1+1} \\ i, k \in V_1} }
	\Delta_{V_1 \cup\{k\}} \left(\kappa^{\geq k+1}_{|B_1|} [f] \right)\prod_{t=2}^p S_{v_t+1}\left(\kappa^{\geq k+1}_{|B_t|} [f] \right)
	\\	&\leq  \sum_{\substack {v_1,\ldots, v_p \in \mathbb{N} \\ v_1 +\cdots +v_p = v}} 
	\sum_{k=j+1}^n  \sum_{\substack{V_1 \in \binom{[n]}{v_1+1} \\ i, k \in V_1} }
	\Delta_{V_1 \cup\{k\}} \left(\kappa^{\geq k+1}_{|B_1|} [f] \right)\prod_{t=2}^p \hbar(v_t+1, |B_t|)
	\\ &\leq \sum_{\substack {v_1,\ldots, v_p \in \mathbb{N} \\ v_1 +\cdots +v_p = v}} 
	\prod_{t=1}^p \hbar(v_t+1, |B_t|). \ 
\end{align*}
Since $p \geq 2$, we have that  $|B_t| \leq r-1$ and
\[
 1 \leq  v_t + 1 \leq    v+1 + (r -1 - | B_t|)  \leq m -b_t +  1.
\]
Therefore, the application of the induction hypothesis above is correct.
Estimating similarly  the contribution of the cases when $i\in V_2, \ldots,  i\in V_p$, we obtain
\begin{align*}
	S_v(\kappa_{r}^{\geq j+1}[f])  &\leq 
	\sum_{p=2}^r
	\sum_{
		\{B_1,\ldots,B_p\}\in \mathcal P_r
	}   \(\dfrac32\)^p   p!
\sum_{\substack {v_1,\ldots, v_p \in \mathbb{N} \\ v_1 +\cdots +v_p = v}} 
\prod_{t=1}^p \hbar(v_t+1, |B_t|)\\
&=
\sum_{p=2}^r
\sum_{\substack {b_1,\ldots, b_p \geq 1 \\ b_1 +\cdots +b_p = r}}   \(\dfrac32\)^p    \binom{r}{b_1,\ldots ,b_p}
\sum_{\substack {v_1,\ldots, v_p \in \mathbb{N} \\ v_1 +\cdots +v_p = v}} 
\prod_{t=1}^p \hbar(v_t+1, b_t)
\\
&=
\sum_{p=2}^r
\sum_{\substack {b_1,\ldots, b_p \geq 1 \\ b_1 +\cdots +b_p = r}}   \(\dfrac32\)^p     
 r!
\sum_{\substack {v_1,\ldots, v_p \in \mathbb{N} \\ v_1 +\cdots +v_p = v}} 
\prod_{t=1}^p  \Dfrac{56^{b_t -1}}{b_t^2} \binom{v_t + 2b_t -2}{b_t -1} \alpha^{b_t}
\\&=
\hbar(v,r) \frac{r^2}{\binom{v+2r-3}{r-1}} \sum_{p=2}^r
  \(\dfrac32\)^p      56^{1-p} 
\sum_{\substack {b_1,\ldots, b_p \geq 1 \\ b_1 +\cdots +b_p = r}} 
\sum_{\substack {v_1,\ldots, v_p \in \mathbb{N} \\ v_1 +\cdots +v_p = v}} 
\prod_{t=1}^p   \Dfrac{1}{b_t^2}\binom{v_t + 2b_t -2}{b_t -1}. 
\end{align*}




Note that $\binom{v_t+2b_t-2}{b_t-1}$ 
is the number of ways to
write $v_t+b_t-1$ as the sum of $b_t$ nonnegative integers.
Some sums don't occur; for example, the sum 1 is impossible for $b_t=3,v_t\ge 0$.
Therefore,
\[
\sum_{\substack {v_1,\ldots, v_p \in \mathbb{N} \\v_1 +\cdots +v_p = v}}
 \prod_{t=1}^p\binom{v_t + 2b_t -2}{b_t -1}
 \le \binom{v + 2r - p - 1}{r-1}
 \le \binom{v + 2r - 3}{r-1},
\]
since the middle expression is the number of ways to write $\sum_t (v_t+b_t-1)=v+r-p$ as the sum of $\sum_t b_t=r$ nonnegative integers,
and the last expression follows since $p\ge 2$. 

Next, by induction on $p\geq 2$, we have that
\begin{equation}\label{eq:sum-b}
	\sum_{\substack {b_1,\ldots, b_p \geq 1 \\ b_1 +\cdots +b_p = r}} 
	\prod_{t=1}^p   \Dfrac{1}{b_t^2} \leq  \Dfrac{(5.3) ^{p-1}}{r^2}.
\end{equation}
Both the base of induction and the induction step for \eqref{eq:sum-b} rely on the following bound:
\[
	\sum_{t=1}^{r-1} \Dfrac{r^2}{t^2(r-t)^2} 
	=  \sum_{t=1}^{r-1}   \left(\Dfrac{1}{t} + \Dfrac{1}{r-t}\right)^2 
	\leq \sum_{t=1}^{r-1} \left(\Dfrac{1}{t^2} + \Dfrac{1}{(r-t)^2} + \Dfrac{2}{r-1}\right)\leq 2+ \Dfrac{\pi^2}{3} <5.3.
\]
Combining the above estimates, we conclude the 
\[
S_v(\kappa_{r}^{\geq j+1}[f])   
\leq  4\hbar(v,r)\sum_{p=2}^r   \(\dfrac32\)^p      \left(\dfrac{5.3}{56}\right)^{p-1}
\leq \hbar(v,r).
\]
Thus, we have established the induction step and proved \eqref{ind:hyp}.

Finally, using \eqref{ind:hyp}, Lemma \ref{L:E-Delta},  and bounding 
$\binom{2r-2}{r-1} \leq 2^{2r-2}$, we obtain
\[
	\|\kappa_r^{\geq j+1}[f] - \Ex^j \kappa_r^{\geq j+1}[f]\|_\infty
	\leq  \Delta_{j}(\kappa_r^{\geq j+1}[f]) \leq \hbar(1,r) \leq 80^{r-1} \frac{(r-1)!}{r} \alpha^r
\]
as claimed.
\end{proof}

 \begin{lemma}\label{l:d-cumulant2}
 		Suppose $f \in 	\calF^{\alpha}_m(\X)$ for some   $\alpha \geq 0$ and positive integer $m$.
 	Then, for any   $r \in [m]$ and  $j\in [n]$,
 	we have 
 		\[ \left\|\kappa_{r}^{\geq j +1 } [f]  -\kappa_{r}^{\geq j}[f]\right\|_\infty \leq   1.1\cdot 80^{r-1} \frac{(r-1)!}{r} \alpha^r.\]
 	\end{lemma}
 \begin{proof} 
 	From Lemma \ref{l:cum_rec}(a), we have that
 	\[
 	{\kappa^{\geq j}_r [f] - \Ex^{j}\kappa_{r}^{\geq j +1 } [f]}
 	=   \sum_{p=2}^r\sum_{(B_1,\ldots,B_p)\in \mathcal P_r} 
 	\kappa^{\geq j} \left[\kappa^{\geq j+1}_{|B_1|}[f],
 	\ldots,
 	\kappa^{\geq j+1}_{|B_p|}[f]\right].
 	\] 
 	Recalling from \eqref{cumulant-red} that
 	\[
 	\kappa^{\geq j} \left[\kappa^{\geq j+1}_{|B_1|}[f],
 	\ldots,
 	\kappa^{\geq j+1}_{|B_p|}[f]\right]
 		=                  
 		\kappa^{\geq j}	\left[ \kappa^{\geq j+1}_{|B_1|}[f] - \Ex^j \kappa^{\geq j+1}_{|B_1|}[f]
 		\ldots,
 		\kappa^{\geq j+1}_{|B_p|}[f] - \Ex^j \kappa^{\geq j+1}_{|B_p|}[f]\right],
 	\]
  applying Lemma \ref{cumulantgeneral}(d) with $V= \emptyset$, and using \eqref{eq:sum-b}, we find that
  \begin{align*}	
  			\|	\kappa^{\geq j}_r [f] &- \Ex^{j}\kappa_{r}^{\geq j +1 } [f]
  			\|_\infty \\
  			&\leq \sum_{p=2}^r 
  			\sum_{\{B_1,\ldots,B_p\} \in \mathcal P_r}   \(\dfrac32\)^p (p-1)! \prod_{t=1}^p \|\kappa^{\geq j+1}_{|B_t|}[f] - \Ex^j \kappa^{\geq j+1}_{|B_t|}[f]\|_\infty
  			\\ 
  			&\leq \sum_{p=2}^r 
  			\sum_{\substack {b_1,\ldots, b_p \geq 1 \\ b_1 +\cdots +b_p = r}}
			  \binom{r}{b_1,\ldots, b_p} \(\dfrac32\)^p  \Dfrac1p \prod_{t=1}^p 80^{b_t-1} \dfrac{(b_t-1)!}{b_t} \alpha^{b_t}
  			\\
  			&= 
  			80^{r-1} \frac{(r-1)! }{r} \alpha^r  \sum_{p=2}^r 
  			r^2  \(\dfrac32\)^p   \,80^{1-p} \Dfrac{1}{p} 
  				\sum_{\substack {b_1,\ldots, b_p \geq 1 \\ b_1 +\cdots +b_p = r}}  \prod_{t=1}^p \dfrac{1}{b_t^2}\\
  				&\leq 
  				80^{r-1} \frac{(r-1)! }{r} \alpha^r  \sum_{p=2}^r 
  				   \(\dfrac32\)^p   \,80^{1-p} \dfrac{1}{2} 
  			(5.3)^{p-1} \leq 0.1 \cdot 80^{r-1} \frac{(r-1)! }{r} \alpha^r,
  \end{align*}
where we used \eqref{eq:sum-b}
in the second last inequality.
The bound on  
  			$\|	\kappa^{\geq j+1}_r [f] - \Ex^{j}\kappa_{r}^{\geq j +1 } [f]
\|_\infty$ from  Lemma \ref{l:d-cumulant1} 
and the triangle inequality complete the proof.
 \end{proof}

 \begin{lemma}\label{l:cumulants}
 		Suppose $f \in 	\calF^{\alpha}_m(\X)$ for some   $\alpha \geq 0$ and positive integer $m$.
 	Then, for any      $j\in [n]$,
 	we have 
 		\[
 		\left\|\Ex^{\geq j}\exp\left(\sum_{r=1}^m 
 		\frac{\kappa_{r}^{\geq j +1 } [f]  -\kappa_{r}^{\geq j }[f]} {r!}
 		\right)-1\right\|_\infty \leq e^{ (100\alpha)^{m+1} }-1.
 		\]
 \end{lemma}
 \begin{proof} 
 	First, if $\alpha\geq 1/100$, then using Lemma \ref{l:d-cumulant2}, we obtain
 	\[
 	\left\|\sum_{r=1}^m
 	\frac{\kappa_{r}^{\geq j +1 } [f]  -\kappa_{r}^{\geq j }[f]} {r!}\right\|_\infty
 	\leq 1.1 \sum_{r=1}^m 80^{r-1} \Dfrac{1}{r^2} \alpha^r
 	\leq 1.1\cdot (\alpha + 80^{m-1} \alpha^{m}) \sum_{r=1}^m \Dfrac{1}{r^2}
 	\leq (100 \alpha)^{m+1}.
 	\]
 	Therefore,
 	\begin{align*}
 			\left\|\Ex^{\geq j}\exp\left(\sum_{r=1}^m 
 		\frac{\kappa_{r}^{\geq j +1 } [f]  -\kappa_{r}^{\geq j }[f]} {r!}
 		\right)-1\right\|_\infty 
 		 &\leq 
 			\left\|\exp\left(\sum_{r=1}^m
 		\frac{\kappa_{r}^{\geq j +1 } [f]  -\kappa_{r}^{\geq j }[f]} {r!}
 		\right)-1\right\|_\infty\\  &\leq \sum_{j = 1}^{\infty} \Dfrac{(100 \alpha)^{j(m+1)}}{j!}
 		= e^{(100 \alpha)^{m+1}} -1.
 	\end{align*}
 	Thus, in the following, we can assume that $\alpha \leq  1/100$.
 	
 	Let 
 	\[
 		F(z):=\Ex^{\geq j}\exp\left(\sum_{r=1}^m 
 		\frac{z^r \(\kappa_{r}^{\geq j +1 } [f]  -\kappa_{r}^{\geq j }[f]\)} {k!} \right) -1.
 	\]
 	and let $f_1, f_2,\ldots \in \calF_\infty(\X)$ denote  the coefficients of its Taylor's expansion:
 	\begin{equation}\label{eq:def_pik}
 		F(z)
 		= \sum_{r=1}^\infty z^r f_r.
 	\end{equation}
 	Due to  Lemma \ref{cumulantgeneral}(b,d), 
  the
 	series  $\sum_{r=1}^{\infty} \frac{z^r\kappa_r^{\geq k}[f]}{r!}$ converges 
 for any $k \in [n]$ and   $z\in \Complexes$ with  
 	$|z| < \frac{2}{3\norm{f}_\infty}$   and
 	\[
 	\exp\left(\sum_{r=1}^{\infty} \frac{z^r\kappa_r^{\geq k}(\pi)}{r!}\right) = \Ex^{\geq k} e^{z f}. 
 	\]
 	Taking $k=j,j+1$, we obtain that
 	\[
 		\Ex^{\geq j} \left[ \exp\left(\sum_{r=1}^{\infty} \frac{z^r \(\kappa_r^{\geq j+1}[f] -\kappa_r^{\geq j}[f]\)}{r!}\right) \right] 
 		=\exp\left(\sum_{r=1}^{\infty} \frac{-z^r  \kappa_r^{\geq j}[f]}{r!}\right)
 		\Ex^{\geq j} \left[\Ex^{\geq j+1} [e^{zf}]\right] =1.
	\]
 	It implies that first $m$ terms in series of \eqref{eq:def_pik} are trivial: $f_1  = \cdots  =f_m \equiv 0$. 
 		Applying Cauchy's integral theorem, we obtain
 	\begin{align*}
 	 \Ex^{\geq j}\exp\left(\sum_{r=1}^m 
 	 \frac{  \(\kappa_{r}^{\geq j +1 } [f]  -\kappa_{r}^{\geq j }[f]\)} {k!} \right) -1 
	 &=
 	 F(1)
 	 =
 	   \sum_{r=m+1}^\infty f_r  \\ &=
 		\frac{1}{2\pi i} \oint  \sum_{k>m+1 }  \frac{1}{z^k}  F(z) dz = \frac{1}{2\pi i}
 		\oint  \frac{F(z)}{(z-1) z^{m+1}}dz,
 	\end{align*}
 	where the integrals are over any contour encircling the origin. We take the circle $\{z\in \Complexes \st |z|=\dfrac{1}{80 \alpha}\}$ as the contour. Using Lemma \ref{l:d-cumulant2} and recalling $\alpha <1/100$, we  observe that for any $z$ on this circle
 	\[
 	  \Dfrac{|F(z)|}{|z-1|} \leq \Dfrac{1.1}{100/80 -1} \sum_{r=1}^m 80^{r-1} \Dfrac{1}{r^2} \alpha^r |z|^r
    <1.
 	\]
 	The required bound follows. 
 \end{proof}

\subsection{Proof of Theorem \ref{cumuThm}}
  	Recalling that $\kappa_1^{\geq n+1} [f] = f$ and $\kappa_r^{\geq n+1}[f]=0$ 
  	   for all $r\geq 2$, we have
	   \[
  	   		f(\X) -  \sum_{r=1}^m \frac{\kappa_r (f(\X))}{r!}
  	   	=  \sum_{j=1}^n  \sum_{r =1}^m\frac{\kappa_{r}^{\geq j +1 } [f]  -\kappa_{r}^{\geq j }[f]} {r!}.
  	   \]
  	   Applying Lemma \ref{l:cumulants}, we find that 
  	   \begin{align*}
  	  &\Ex \left[ \exp \left( \sum_{j=1}^n  \sum_{r =1}^m \frac{\kappa_{r}^{\geq j +1 } [f]  -\kappa_{r}^{\geq j }[f]} {r!} \right) \right]\\
  	  & {\quad}=
  	   \Ex \left[  \exp \left( \sum_{j=1}^{n-1} \sum_{r=1}^m \frac{\kappa_{r}^{\geq j +1 } [f]  -\kappa_{r}^{\geq j }[f]} {r!} \right) \cdot \Ex^{\geq n} \left[
  	   \exp\left(\sum_{r=1}^m \frac{\kappa_{r}^{\geq n +1 } [f]  -\kappa_{r}^{\geq n }[f]} {r!}\right) \right]\right]
  	   \\
  	   & {\quad}=   (1+\delta_1)\Ex \left[  \exp \left( \sum_{j=1}^{n-1} \sum_{r=1}^m \frac{\kappa_{r}^{\geq j +1 } [f]  -\kappa_{r}^{\geq j }[f]} {r!} \right) \right]= \cdots = \prod_{j=1}^n \,(1+\delta_j),
  	   \end{align*}
  	   where $|\delta_j| \leq e^{(100\alpha)^{m+1}} -1$ for all $j \in [n]$.
      Since $f$ is real-valued,  we also have that $\delta_j>-1$. 
  	  Therefore, 
  	  \[
  	  		\Ex \left[\exp \left(f(\X) -  \sum_{r=1}^m \frac{\kappa_r (f(\X))}{r!} \right)\right]
  	  		= (1+\delta)^n,
  	  \]
  	  where $\delta :=  \left(\prod_{j=1}^n (1+\delta_j)\right)^{1/n} -1$.  Observing that  $\min_{j\in [n]} \delta_j \leq \delta   \leq  \max_{j\in [n]} \delta_j$, we have established that  $|\delta| \leq e^{(100\alpha)^{m+1}} -1$ as required.

 By Lemma \ref{l:d-cumulant2},
we have that
for any $r \in [m]$ and  $j\in [n]$,
\[
\left\|\kappa_{r}^{\geq j +1 } [  f(\X) ]  
- \kappa_{r}^{\geq j}[  f(\X) ] \right\|_\infty 
\leq   1.1\cdot 80^{r-1} \Dfrac{(r-1)!}{r} \alpha^r
\leq 
\Dfrac{(r-1)!}{50r} (80\alpha)^r.
\]
Using the triangle inequality,
\[
| \ka_{r}\left( f(\X) \right) |
\le 
\sum_{j=0}^{n-1}\, \left\|\kappa_{r}^{\geq j +1 } [  f(\x) ]  
- \kappa_{r}^{\geq j}[  f(\x) ] \right\|_\infty, 
\]
and the claimed bound follows.
 
\nicebreak

\nicebreak
\section*{Appendix}

Exact values of $\RT(n)$ for $n\le 21$ were found by McKay~\cite{mckay1983applications} using
a method of summing over roots of unity.  We have extended this to $n\le 37$ using the following
recurrence.
Let $T(d_1,\ldots,d_n)$ denote the number of tournaments on $n$ vertices
with out-degrees $d_1,\ldots,d_n$,
where a tournament is an orientation of an undirected complete graph.
Then, by considering the removal of vertex $n$,
\[
   T(d_1,\ldots,d_n) = \sum_{\substack{e_1,\ldots,e_{n-1}\in\{0,1\}\\ e_1+\cdots+e_{n-1}=n-1-d_n}}
     \kern-1em T(d_1-e_1,\ldots,d_{n-1}-e_{n-1}),
\]
and also note that $T$ is a symmetric function of its arguments.
See Table~\ref{Table:RT}.

\begin{table}[ht]
\small
\centering
\begin{tabular}{r|p{0.87\textwidth}}
$n$ & $\RT(n)$ \\
\hline
1 & 1 \\
3 & 2 \\
5 & 24 \\
7 & 2640 \\
9 & 32 30080 \\
11 & 4 82515 08480 \\
13 & 9 30770 06112 92160 \\
15 & 240 61983 49824 94283 79648 \\
17 & 85584 72055 41481 49511 79758 79680 \\
19 & 4271 02683 12628 45202 01657 80015 93666 76480 \\
21 & 3035 99177 67255 01434 06909 90026 40396 04333 20198 14400 \\
23 & 31111 25335 58482 03432 16879 55029 99798 94772 74014 27415 01378 56000 \\
25 & 46 41175 34102 33590 76153 19841 21486 62289 71154 35036 87620 35567 90979
81774 06976 \\
27 & 1 01613 79494 93595 16286 17799 57707 58654 34480 25582 38882 79881 93794
83077 47797 07683 48315 64800 \\
29 & 32874 42487 57440 71437 03099 54561 70730 45735 81860 90338 16899 44155 
56007 47117 93383 29715 53931 99301 61725 44000 \\
31 & 1 58080 31329 37879 48113 48365 61225 94846 34453 23284 20116 08717 95271
39910 78379 57102 46216 62189 55225 10623 95890 69751 74000 64000 \\
33 & 113 55331 66724 13409 50706 27943 25155 74835 60333 94267 70339 22296
31722 43704 16674 54847 32916 31865 21307 04028 02521 10357 42226 17221
58828 67015 68000 \\
35 & 1 22386 44546 20140 32917 57098 60021 24725 82639 58465 65811 80615
89271 23479 00218 68020 11957 16753 08650 67421 64775 78535 34744 86223
90705 83144 23519 70257 72801 31373 46560 \\
37 & 19868 18615 30379 61435 68362 02070 63930 19820 07449 69481 15232
49050 85973 12501 88142 39611 14080 38454 82389 53944 91763 75862 81845
56861 43415 69406 38026 42548 47415 17715 34651 86485 95415 04000
\end{tabular}
\normalsize
\caption{Counts of labelled regular tournaments\label{Table:RT}}
\end{table}

\end{document}